\documentclass[11pt]{amsart}
\usepackage{amsmath,amssymb,amsthm}

\newtheorem{theorem}{Theorem}
\newtheorem{lemma}[theorem]{Lemma}
\newtheorem{corollary}[theorem]{Corollary}

\newcommand{\SSfootnote}[1]{\renewcommand{\baselinestretch}{1}%
             \footnote{#1}\renewcommand{\baselinestretch}{1}}

\newcommand{\zbar}{\bar{z}}
\newcommand{\wbar}{\bar{w}}
\newcommand{\del}[1]{\frac{\partial}{\partial #1}}
\newcommand{\delsq}[1]{\frac{\partial^{2}}{\partial #1^{2}}}
\newcommand{\x}{\chi}

\newcommand{\xbar}{\bar{x}}
\newcommand{\half}{\frac{1}{2}}

\renewcommand{\Im}{\mbox{\rm Im\,}}
\renewcommand{\Re}{\mbox{\rm Re\,}}

\newcommand{\R}{\mathbb{R}}
\renewcommand{\H}{\mathbb{H}}
\newcommand{\C}{\mathbb{C}}

\newcommand{\p}{\Phi}
\newcommand{\dbar}{\bar{\partial}}

\renewcommand{\r}{\rho}

\newcommand{\SO}{\mbox{\rm SO}}
\newcommand{\SU}{\mbox{\rm SU}}
\newcommand{\Sp}{\mbox{\rm Sp}}
\newcommand{\U}{\mbox{\rm U}}
\newcommand{\Spin}{\mbox{\rm Spin}}

\newcommand{\so}{\mathfrak{so}}
\renewcommand{\sp}{\mathfrak{sp}}
\newcommand{\spin}{\mathfrak{spin}}
\newcommand{\su}{\mathfrak{su}}
\renewcommand{\u}{\mathfrak{u}}

\begin{document}

\title{Singular Instantons with $\SO(3)$ Symmetry}
\author{Gregory D. Landweber}

\email{greg@math.uoregon.edu}
\urladdr{http://www.uoregon.edu/\~{}greg/}
\address{Mathematics Department\\
         University of Oregon\\
         Eugene, OR 97403-1222}

\begin{abstract}
	The purpose of this article is to provide an explicit construction
	for a family of singular instantons on $S^{4} \setminus S^{2}$ with
	arbitrary real holonomy parameter $\alpha$. This family includes
	the original $\alpha = 1/4$, $c_{2} = 3/2$ solution discovered by
	P.~Forg\'{a}cs, Z.~Horv\'{a}th, and L.~Palla, and our approach is
	modeled on that of their 1981 paper. Our primary tool is the ansatz
	due to Corrigan, Fairlie, Wilczek, and 't Hooft that constructs a
	self-dual Yang-Mills connection using a positive real-valued harmonic
	super-potential. Here we reformulate this harmonic function ansatz
	in terms of quaternionic notation, and we show that it arises naturally
	from the Levi-Civita connection of a conformally Euclidean metric.
	
	To simplify the construction, we introduce an $\SO(3)$-action on
	$S^{4}$, and we show by dimensional reduction that the symmetric
	self-duality equation on $S^{4}$ is equivalent to the vortex equations
	over hyperbolic space $\mathcal{H}^{2}$. We thus obtain a similar
	harmonic function ansatz for hyperbolic vortices, which we also
	derive using conformal transformations of $\mathcal{H}^{2}$. Using
	this ansatz, we construct the vortex equivalents of the symmetric
	't Hooft instantons, and we prove using the equivariant ADHM
	construction that they provide a complete description of all
	hyperbolic vortices. We also analyze when two vortices constructed
	by this ansatz are gauge equivalent, obtaining the surprising result
	that two such vortices are completely determined by the gauge 
	transformation between them.
\end{abstract}

\maketitle


\section*{Introduction}

In recent years, the study of singular Yang-Mills fields has been an 
extremely active area of research. Considering $\SU(2)$ instantons on
four manifolds with codimension two singularities, it was found that
these connections can admit non-trivial holonomy around arbitrarily
small circles linking the embedded singular surface. An analytical theory
for such instantons with holonomy singularity has been developed in
\cite{KM,KM2}. Although we currently have an understanding
of the moduli space for such singular instantons, the literature in this
field has a conspicuous dearth of explicit examples. A singular solution on
$T^{2}\times D^{2}$ is given in the appendix to \cite{KM}, although of 
greater interest are solutions on the standard model $S^{4}\setminus S^{2}$.

The first example of a singular instanton was discovered by P.~Forg\'{a}cs,
Z.~Horv\'{a}th, and L.~Palla.  In their 1981 paper \cite{FHP1}, they 
describe a self-dual Yang-Mills field on $S^{4}\setminus S^{2}$ with 
the fractional Chern class $c_{2} = 3/2$.  At first their result was not 
readily accepted, due to the resistance to the new idea of a 
fractional charge.
They later published a second paper
\cite{FHP2} defending their result, and since then over a decade of successful
research into the field has eliminated any initial skepticism. 
Nevertheless, their construction itself remains poorly understood.

The goal of this article is to elucidate and extend the work of
Forg\'{a}cs {\em et al.\/}, writing their construction using simpler 
notation, explaining the motivation behind their 
formulae,
and generalizing to obtain a family of singular instantons with varying
holonomy parameter.
To this work we shall contribute a
mathematical perspective,
exchanging indices and Pauli matrices for more invariant complex, 
quaternionic, and spinor notation, and offering geometric interpretations
for the equations involved.

Section 1 is devoted to the construction of instantons on $S^{4}$ employing
the ansatz proposed by the physicists Corrigan, Fairlie, and Wilczek in
1976 and described in \cite{CF, JNR}. Starting with a positive real-valued 
function $\r$ on $\R^{4}$, known as the {\em super-potential}, we 
consider the Yang-Mills connection
\begin{equation} \label{eq:ansatz}
	A = \sum_{\mu,\nu} i\bar{\sigma}_{\mu\nu}\,\partial^{\nu}\log\r\,dx^{\mu}.
\end{equation}
Here the anti-symmetric matrix $\bar{\sigma}_{\mu\nu}$ is defined as%
\SSfootnote{We use the Greek indices $\mu,\nu$ when indexing over 4-space,
            while the Roman indices $i,j$ range over 1,2,3.}
\begin{displaymath}
	\bar{\sigma}_{\mu\nu} = \left\{ \begin{array}{l}
	                          \bar{\sigma_{ij}} = \frac{1}{4i}
	                              \left[ \sigma_{i},\sigma_j \right] \\
	                          \bar{\sigma_{i0}} =
	                              -\frac{1}{2}\sigma_{i}
	                          \end{array}
	                  \right.
\end{displaymath}
where the $\sigma_{i}$ are the standard Pauli matrices generating the Lie
algebra $\su(2)$. For such a connection, the self-duality equation 
$\ast F_{A} = F_{A}$ is equivalent to the condition $\Delta\rho = 0$.
By reversing orientation, this construction can also be used to generate
anti-self-dual connections from a harmonic super-potential.

This harmonic function ansatz was used by 't Hooft to construct a class
of instantons with $5n$ parameters, corresponding to the centers and
scales of $n$ superimposed basic instantons. Since then, this ansatz
has been shown to be the simplest case of a more general
algebraic-geometric construction involving twistors discussed in \cite{AW}.
More recently, both constructions have been eclipsed by the ADHM
description of instantons given in \cite{ADHM}, which provides a complete construction for all ASD connections on $S^{4}$ up to gauge equivalence.
In Section 1 we recast the harmonic function ansatz in terms of 
quaternionic notation. Not only does this greatly simplify the
required calculations, but also it better exhibits the underlying
structure. We then show how these connections arise naturally via
conformal transformations.

In Section 2 we introduce an $\SO(3)$-action on $S^{4}$. Taking
advantage of the conformal equivalence
$S^{4}\setminus S^{1} \cong \mathcal{H}^{2}\times S^{2}$, we show that the
symmetric SD and ASD equations over $S^{4}$ are equivalent to the vortex
and anti-vortex equations over hyperbolic space $\mathcal{H}^{2}$. This
technique is known as dimensional reduction. The harmonic function
ansatz for instantons then reduces to a similar ansatz for hyperbolic
vortices, which we also derive using conformal transformations of 
hyperbolic space. After computing the vortex equivalents of the
symmetric 't Hooft instantons, we use an equivariant version
of the ADHM construction to provide a classification for all hyperbolic
vortices. Examining gauge transformations, we obtain the surprising
result that if two hyperbolic vortices constructed by the harmonic
function ansatz are gauge equivalent, then they are both completely
determined by the gauge transformation between them.

We return to our primary task of constructing singular instantons in 
Section 3. Restricting our attention to $\SO(3)$-invariant connections
on $S^{4}$, we can work instead with hyperbolic vortices. Using the unit 
disc model of $\mathcal{H}^{2}$, singular instantons correspond to vortices
with a holonomy singularity at the origin. We then proceed to construct
solutions on the cut disc using the harmonic function ansatz, patching them 
together with gauge transformations to form global solutions on the 
punctured disc. In \S\ref{fhp} we essentially rewrite the paper 
\cite{FHP1} in this context, and in the following section we construct
our desired family of singular vortices. 

\section{The Harmonic Function ansatz}

\label{instanton-ansatz}

\subsection{Quaternionic Notation}

\label{quaternionic-notation}

For the duration of this section, we adopt the quaternionic notation as used
in \cite{A}. Writing $x\in\H$ in the form $x = x^{0}+ix^{1}+jx^{2}+kx^{3}$,
its conjugate is $\xbar = x^{0}-ix^{1}-jx^{2}-kx^{3}$, and the
corresponding differentials are
\begin{align*}
    dx      =  dx^{0} + i\,dx^{1} + j\,dx^{2} + k\,dx^{3} \qquad
    d\xbar  =  dx^{0} - i\,dx^{1} - j\,dx^{2} - k\,dx^{3}.
\end{align*}
By analogy with the complex case, we define the partial derivatives
\begin{align*}
    \del{x}     & = \half \left(
                  \del{x^{0}} - i\del{x^{1}} - j\del{x^{2}} - k\del{x^{3}}  
                            \right) \\
    \del{\xbar} & = \half \left(
                  \del{x^{0}} + i\del{x^{1}} + j\del{x^{2}} + k\del{x^{3}}  
                            \right).
\end{align*}
In this notation the Laplacian takes the form
\begin{displaymath}
	\Delta = -4\,\del{x}\del{\xbar} = -4\,\del{\xbar}\del{x}.
\end{displaymath}
As expected, the exterior derivative $d$ may be written as the sum of
$\partial$ and $\dbar$ components, although there are now two distinct
splittings
\begin{displaymath}
	d =  dx\,\del{x} + \del{\xbar}\,d\xbar 
	  =  \del{x}\,dx + d\xbar\,\del{\xbar} 
\end{displaymath}
due to the non-abelian nature of the operators involved. Expanding the
2-form $dx\wedge d\xbar$ in terms of coordinates as
\begin{equation}\begin{split}
	dx\wedge d\xbar & =  -2 \left[\,
	    i\left(dx^{0}\wedge dx^{1} + dx^{2}\wedge dx^{3}\right)\,+\,
	    j\left(dx^{0}\wedge dx^{2} + dx^{3}\wedge dx^{1}\right)
	                         \right. \nonumber \\
	&  \qquad  \left.\mbox{\qquad}\,+\,             
	    k\left(dx^{0}\wedge dx^{3} + dx^{1}\wedge dx^{2}\right)
	                     \,\right], \label{eq:dxdxbar}
\end{split}\end{equation}
we see that $dx\wedge d\xbar$ is self-dual and likewise that
$d\xbar\wedge dx$ is anti-self-dual. 
 
Rewriting the connection (\ref{eq:ansatz}) in terms of this new 
quaternionic notation, the harmonic function ansatz now takes the
surprisingly familiar form

\begin{theorem} \label{theorem-1}
    Given a positive real-valued super-potential $\r$ on $\R^{4}$,
    the Yang-Mills connection $A^{+}$ defined by
	\begin{equation}
	    \label{eq:anti-self-dual}
	    A^{+}= - \Im\left(\del{\xbar}\log \r\,d\xbar \right)
	         = - \half \left( \del{\xbar}\log \r\,d\xbar -
	                           dx\,\del{x}\log\r \right)
	\end{equation}
	is anti-self-dual and the connection $A^{-}$ defined by the conjugate
	expression
	\begin{equation}
	   \label{eq:self-dual}
	    A^{-} = - \Im \left(\del{x}\log \r\,dx \right)
	          = - \half \left( \del{x}\log \r\,dx -
	                           d\xbar\,\del{\xbar}\log\r \right)
	\end{equation}
	is self-dual if and only if the super-potential $\r$ is harmonic.
\end{theorem}

Before proceeding with the proof of this theorem the reader may want
to verify that (\ref{eq:ansatz}) and (\ref{eq:self-dual}) both yield the
same self-dual connection.  Expanding (\ref{eq:self-dual}) using
coordinates, we obtain the expression
\begin{equation*}\begin{split}
	A^{-} & =  \half \left(\,
	        (+i\partial_{1}+j\partial_{2}+k\partial_{3})\,dx^{0} +
	        (-i\partial_{0}-k\partial_{2}+j\partial_{3})\,dx^{1}
	              \right. \nonumber\\
	      &  \qquad       \left. \,\mbox{}+
	        (-j\partial_{0}+k\partial_{1}-i\partial_{3})\,dx^{2} +
	        (-k\partial_{0}-j\partial_{1}+i\partial_{2})\,dx^{3}
	              \,\right),
\end{split}\end{equation*}
writing $\partial_i$ as an abbreviation for $\partial_i\log\r$.

\begin{proof}[Proof of Theorem 1]
For the purposes of this proof, we restrict
our attention to the potentially self-dual connection $A^{-}$ given in
(\ref{eq:self-dual}), calling it $A$. Explicitly computing the two
components of the curvature $F_{A} = dA+ A\wedge A$, we obtain
\begin{align*}
	A\wedge A & = -\half \left(
	                 \del{x}\log\r\,dx \wedge d\xbar\,\del{\xbar}\log\r +
	                 d\xbar\,\del{\xbar}\log\r \wedge \del{x}\log\r\,dx
	                       \right) \\
	dA & =  -\half \left( \del{x}\,dx + d\xbar\,\del{\xbar} \right)
	                \left(
	                      \del{x}\log \r\,dx - d\xbar\,\del{\xbar}\log\r
	                \right) \\
	   & =  -\half \left(
	                - \del{x}\,dx\wedge d\xbar\,\del{\xbar}\log\r +
	                  d\xbar\,\del{\xbar}\wedge\del{x}\log\r\,dx
	                \right) .
\end{align*}
Recalling that the 2-forms $dx\wedge d\xbar$ and $d\xbar\wedge dx$ are 
self-dual and anti-self-dual respectively, the curvature of $A$ then splits
as $F_{A} = F_{A}^{+} + F_{A}^{-}$ with
\begin{align*}
    F_{A}^{+} & =  \half \left(
                            \del{x}\,dx\wedge d\xbar\,\del{\xbar}\log\r\,-\,
                            \del{x}\log\r\,dx \wedge d\xbar\,\del{\xbar}\log\r
                          \right) \\
    F_{A}^{-} & =  -\half \left(
                            \del{\xbar}\del{x}\log\r\,+\, 
                            \del{\xbar}\log\r\,\del{x}\log\r
                           \right) d\xbar\wedge dx.
\end{align*}
Noting the identity
\begin{displaymath}
	\del{\xbar}\del{x}\log\r + \del{\xbar}\log\r\,\del{x}\log\r
	  = \frac{1}{\r}\,\del{\xbar}\del{x}\r = -\frac{1}{4\r}\Delta\r,
\end{displaymath}
we see that the self-duality equation $F_{A}^{-} = 0$ is equivalent to
the condition $\Delta\r = 0$ that the super-potential $\r$ be harmonic.
\end{proof}


We now calculate the curvature density $|F_{A}|^{2}$ of the self-dual
connection (\ref{eq:self-dual}), from which we can construct the Yang-Mills
functional $\|F_{A}\|^{2}$ and the Chern class $c_{2}(A)$. Using the
decomposition $F_{A} = F_{A}^{+} + F_{A}^{-}$ given above, we first 
compute the anti-self-dual component $|F_{A}^{-}|^{2}$.  From equation
(\ref{eq:dxdxbar}) we observe that
$(d\xbar\wedge dx)\wedge-(\overline{d\xbar\wedge dx}) = 24\,d\mu$,
where $d\mu$ is the volume form, and we immediately obtain
\begin{equation*}
	|F_{A}^{-}|^{2} = \frac{3}{8}\left(\frac{1}{\r}\Delta\r\right)^{2},
\end{equation*}
which clearly vanishes if the super-potential $\r$ is harmonic.


\newcommand{\logr}{\log\r}
\newcommand{\di}{\partial_{i}}
\renewcommand{\dj}{\partial_{j}}

On the other hand, the self-dual component $|F_{A}^{-}|^{2}$ of the 
curvature density is significantly more difficult to compute. Again
using the expansion (\ref{eq:dxdxbar}) for $dx\wedge d\xbar$, we have
\begin{equation*}\begin{split}
	|F_{A}^{+}|^{2} & =  2 \left(
	      \left|
	    \del{x}\,i\,\del{\xbar}\log\r - \left(\del{x}\log\r\right)
	                                   i\left(\del{\xbar}\log\r\right)
	      \right|^{2}\right. \\
	&   \qquad \mbox{}+\left.\left|
	    \del{x}\,j\,\del{\xbar}\log\r - \left(\del{x}\log\r\right)
	                                   j\left(\del{\xbar}\log\r\right)
	     \right|^{2}\right. \\
	&   \qquad \mbox{}+\left.\left|
	    \del{x}\,k\,\del{\xbar}\log\r - \left(\del{x}\log\r\right)
	                                   k\left(\del{\xbar}\log\r\right)
	      \right|^{2}\right),
\end{split}\end{equation*}
which when fully expanded in terms of coordinates becomes
\begin{equation*}\begin{split}
    |F_{A}^{+}|^{2}
	& =  \frac{1}{8}\,\sum_{i,j}\left[
	          4\,(\di\dj\logr)^{2} +
	          3\,(\di\logr)^{2}(\dj\logr)^{2} \right. \\*
	&    \mbox{\qquad} - \left.
	          8\,(\di\dj\logr)(\di\logr)(\dj\logr) -
	          (\di^{2}\logr)(\dj^{2}\logr) \right. \\
	&    \rule{0in}{4ex}\mbox{\qquad} + \left.
	          (\di^{2}\logr)(\dj\logr)^{2} +
	          (\di\logr)^{2}(\dj^{2}\logr)
	                 \right].
\end{split}\end{equation*}
If the super-potential $\r$ is harmonic, then we can take advantage
of the identity $\sum_{i}\di^{2}\log\r = -\sum_{i}(\di\log\r)^{2}$
to simply this expression for $|F_{A}^{+}|^{2}$ to
\begin{displaymath}
	|F_{A}^{+}|^{2} = \frac{1}{2}\,\sum_{i,j}\left[
	          (\di\dj\logr)^{2} -
	          2\,(\di\dj\logr)(\di\logr)(\dj\logr)
	                             \right].
\end{displaymath}
On the other hand, expanding the expression $\Delta\Delta\log\r$, we obtain
\begin{equation*}\begin{split}
	\Delta\Delta\log\r
	& =  \sum_{i,j}\dj^{2}\di^{2}\logr
	= -\sum_{i,j}\dj^{2}(\di\logr)^{2} \\
	& =  -2\,\sum_{i,j}\dj\left[(\di\logr)(\di\dj\logr)\right] \\ 
	& =  -2\,\sum_{i,j}\left[
	      (\di\dj\logr)^{2} + (\di\logr)(\di\dj^{2}\logr)
	                  \right] \\
	& =  -2\,\sum_{i,j}\left[
	      (\di\dj\logr)^{2} - (\di\logr)\di(\dj\logr)^{2}
	                  \right] \\
	& =  -2\,\sum_{i,j}\left[
	      (\di\dj\logr)^{2} - 2\,(\di\logr)(\dj\logr)(\di\dj\logr)
	                  \right],
\end{split}\end{equation*}
again assuming that $\r$ is harmonic and using the 
identity $\sum_{i}\di^{2}\log\r = -\sum_{i}(\di\log\r)^{2}$ repeatedly.


Hence if the super-potential $\r$ is harmonic, then the components of
the curvature density $|F_{A}|^{2}$ for the self-dual connection
(\ref{eq:self-dual}) are
\begin{displaymath}
    |F_{A}^{+}|^{2} = -\frac{1}{4}\,\Delta\Delta\log\r, \qquad
    |F_{A}^{-}|^{2} = 0,
\end{displaymath}
and the Chern class $c_{2}(A)$ and $L^{2}$ norm $\|F_{A}\|^{2}$ are given by
\begin{equation} \label{eq:c2}
	c_{2}(A) = \frac{1}{4\pi^{2}}\,\|F_{A}\|^{2}
	         = -\frac{1}{16\pi^{2}} \int_{\R^{4}}\Delta\Delta\log\r\,d\mu.
\end{equation}
Here we have a factor of $4\pi^{2}$ instead of the customary $8\pi^{2}$
because the function $\xi\mapsto\mbox{Tr}(\xi^{2})$ on the Lie algebras
$\su(2)$ and $\so(3)$ corresponds to the map $\xi\mapsto 2\xi^{2}$
in our quaternionic notation. Similarly, if we take the anti-self-dual
connection (\ref{eq:anti-self-dual}) then the two components
$|F_{A}^{+}|^{2}$ and $|F_{A}^{-}|^{2}$ are interchanged and the Chern
class $c_{2}(A)$ switches sign. It is important to note that the scalar
curvature density is gauge invariant. In other words, if two harmonic
super-potentials $\r_{1}$ and $\r_{2}$ yield gauge equivalent connections
via the ansatz of Theorem~\ref{theorem-1}, then they must satisfy the
equation $\Delta\Delta\log\r_{1} = \Delta\Delta\log\r_{2}$.


\subsection{The 't Hooft Construction}

\label{tHooft}

As an example of the harmonic function ansatz, we take for our
super-potential the Green's functions of the Laplacian.  Although
these functions have $O(1/r^{2})$ poles, the corresponding
singularities can be removed from the resulting connections by a
gauge transformation. In the simplest case, consider the spherically
symmetric harmonic function
\begin{displaymath}
	\r = 1 + \frac{1}{|x|^{2}} = 1 + \frac{1}{x\xbar},
\end{displaymath}
the sum of the Green's functions centered at the origin and infinity.
Applying formula~(\ref{eq:anti-self-dual}), this super-potential
generates the anti-self-dual connection
\begin{equation} \label{eq:singular-gauge}
	A =  \Im\left( \frac{\xbar^{-1}\,d\xbar}{1 + x\xbar} \right)
	  = -\Im\left( \frac{dx\,x^{-1}}{1 + x\xbar} \right),
\end{equation}
which is simply the basic instanton with $c_{2} = 1$ expressed in the 
``singular gauge''.  Applying the gauge transformation $g = x^{-1}$,
we can remove the $O(1/r)$ pole at the origin to obtain this 
instanton's customary form
\begin{equation} \label{eq:standard-gauge}
	g(A) = \Im\left(
	       - x^{-1}\,\frac{dx\,x^{-1}}{1+x\xbar}\,x - dx^{-1}\,x
	          \right)
	     = \Im\left( \frac{\xbar\,dx}{1+x\xbar} \right).
\end{equation}
Note that if we switch to coordinates around infinity by putting
$x = y^{-1}$, then we simply interchange these two gauges
(\ref{eq:singular-gauge}) and (\ref{eq:standard-gauge}). This connection 
therefore takes the same form about infinity as it does about the origin.

More generally, we can modify the basic
instanton~(\ref{eq:standard-gauge}) by applying a dilation
and translation $x\mapsto\,\lambda^{-1}(x-a)$ with $\lambda>0$ real
and $a\in\H$. The super-potential and associated connection then become
\begin{equation}\label{eq:basic-instanton}
	\r = 1 + \frac{\lambda^{2}}{|x-a|^{2}},\qquad
	g_{a}(A) = \Im\left( \frac{ (\xbar - \bar{a})\,dx }{\lambda^{2} + |x-a|^{2}}
	          \right).
\end{equation}
Here we have again used a gauge transformation $g_{a} = (x - a)^{-1}$
in order to remove the singularity at the point $x = a$.


One of the interesting features of this ansatz is that it allows us to take
the superposition of several such instantons simply by adding their 
super-potentials.  For instance, the 't Hooft instantons with $c_{2} = k$
are constructed using the harmonic function
\begin{equation}\label{eq:tHooft}
	\r = 1 + \frac{\lambda_{1}^{2}}{|x-a_{1}|^{2}} + \cdots
	       + \frac{\lambda_{k}^{2}}{|x-a_{k}|^{2}},
\end{equation}
combining $k$ basic instantons of the form (\ref{eq:basic-instanton})
with scales $\lambda_{1},\ldots,\lambda_{k}$ and distinct centers
$a_{1},\ldots,a_{k}$.


\subsection{Conformal Transformations}

\label{conformal-instanton}

In this section, we discuss a differential geometric interpretation of
the harmonic function ansatz introduced in \S\ref{quaternionic-notation}.
Treating the super-potential as a conformal transformation of flat Euclidean 
space, the connections (\ref{eq:anti-self-dual}) and (\ref{eq:self-dual})
arise naturally from the action of the Levi-Civita connection on the
half-spin bundles. We can then express the curvatures of these two 
connections in terms of the decomposition of the Riemann curvature into
its scalar, trace-free Ricci, and conformally invariant Weyl curvature
components, thereby providing an alternative proof of Theorem \ref{theorem-1}.

Starting with the flat Euclidean metric $g_{ij} = \delta_{ij}$ on
$\R^{4}$, we consider the conformally equivalent metric $g' = \r^{2}g$,
given a smooth, positive, real-valued super-potential $\r$. The condition
that $\r$ be harmonic enters when calculating the scalar curvature of this
new metric as in the following lemma.

\begin{lemma} \label{lemma-scalar-curvature}
    The scalar curvature $R'$ of the conformally Euclidean metric $g'$
    given by $g' = \r^{2}\delta_{ij}$ vanishes if and only if the
    super-potential $\r$ is harmonic.
\end{lemma}

\begin{proof}
Using the expression for $R'$ computed in \cite[p. 125]{Au},
in dimension $n=4$ we have
\begin{equation*}\begin{split}
	R'&=- \r^{-2} \left[
	        (n-1)\,\sum_{\nu}\partial_{\nu}^{2}\log\r^{2}\,+\,
	        \frac{(n-1)(n-2)}{4}\,
	        \sum_{\nu}\left(\partial_{\nu}\log\r^{2}\right)^{2}
	               \right] \\
	 &=- 6 \r^{-2}\,\sum_{\nu} \left[ \partial_{\nu}^{2}\log\r +
	                                  \left(\partial_{\nu}\log\r\right)^{2}
	                           \right]
	 =  6 \r^{-3}\Delta\r.\rule{0in}{2.5ex}
\end{split}\end{equation*}
Hence $R' = 0$ if and only if $\Delta\r = 0$.
\end{proof}

Let $\{e_{0},\ldots,e_{3}\}$ be an orthonormal tangent frame for the
original metric $g$. After applying the conformal transformation,
the Levi-Civita connection for the metric $g'$ is given with respect
to this frame by Christoffel's formula
\begin{displaymath}
	\Gamma'^{j}_{ik} = \partial_{i}\log\r\:\delta^{j}_{k} +
	                   \partial_{k}\log\r\:\delta^{j}_{i} -
	                   \partial^{j}\log\r\:\delta_{ik}.
\end{displaymath}
In order to express this as an $\so(4)$ connection, we must rescale the
tangent frame so that it is again orthonormal with respect to the new
metric $g'$. Switching to the frame $e'_{i} = \r^{-1}e_{i}$ introduces
a factor of $-\partial_{i}\log\r\;\delta^{j}_{k}$ into the connection,
cancelling the diagonal term and leaving us with an expression
skew-symmetric in the indices $j$ and $k$.

Taking the double cover $\Spin(4)$ of $\SO(4)$, we recall that the 
Lie algebra isomorphism $\so(4) \cong \spin(4)$ associates to a
skew-symmetric matrix $a_{ij}$ the Clifford algebra element%
\SSfootnote{The Clifford algebra $\mbox{Cl}(4)$ is the algebra generated
            by $\R^{4}$ subject to the relation $v\cdot w + w\cdot v = 
            -2(v,w)$, and the Lie algebra $\spin(4)$ is the subspace
            spanned by $\{e_{i}\cdot e_{j}\}_{i\neq j}$.}
$-\frac{1}{4}\sum_{i,j} a_{ij}\,e_{i}\cdot e_{j}$ (see \cite{LM}).
We may thus write the Levi-Civita connection in this $\spin(4)$
notation as
\begin{displaymath}
	A = \half \sum_{i\neq j}\left( e'_{j}\,\partial^{j}\log\r
	            \right) \cdot \left( e'_{i}\,dx^{i} \right).
\end{displaymath}
From the decomposition $\Spin(4) = \Sp(1) \times \Sp(1)$, we see
that the complex 4-dimensional spin space splits as the direct sum
$S = S^{+} \oplus S^{-}$ of two half-spin spaces, each of which is
isomorphic to the quaternions $\H$. These spaces $S^{+}$ and $S^{-}$
are called the spaces of self-dual and anti-self-dual spinors
respectively. The two half-spin representations $\gamma^{\pm}$ of the
Lie algebra $\spin(4)$ on $\H^{\pm}$ are then given by%
\SSfootnote{Here the various signs are determined by the Clifford algebra
            relation $\gamma^{\pm}(v\cdot w + w\cdot v) = -2(v,w)$
            and also by the convention that 
            $\gamma^{+}(e'_{0}\cdot e'_{1} - e'_{2}\cdot e'_{3}) = 0$
            and
            $\gamma^{-}(e'_{0}\cdot e'_{1} + e'_{2}\cdot e'_{3}) = 0$.}
\begin{align*}
	\gamma^{+}:v\cdot w \mapsto -\gamma(v)\,\gamma^{\ast}(w) \qquad
	\gamma^{-}:v\cdot w \mapsto -\gamma^{\ast}(v)\,\gamma(w),
\end{align*}
where the Clifford action $\gamma(\cdot)$ is simply quaternion 
multiplication
\begin{displaymath}
		\gamma(e'_{0}) = 1, \qquad \gamma(e'_{1}) = i, \qquad
		\gamma(e'_{2}) = j, \qquad \gamma(e'_{3}) = k,
\end{displaymath}
and $\gamma^{\ast}(\cdot)$ is its adjoint. Hence, the Levi-Civita
connection for the conformally transformed metric $g' = \r^{2}g$
splits into the two $\sp(1)$ components
\begin{equation*}
	A^{+}=-\Im\left(\del{\xbar}\log\r\,d\xbar\right) \qquad
	A^{-}=-\Im\left(\del{x}\log\r\,dx\right)
\end{equation*}
acting on the positive and negative half-spin spaces respectively.
Note that these two connections agree with the connections $A^{+}$
and $A^{-}$ given in equations (\ref{eq:anti-self-dual}) and
(\ref{eq:self-dual}).


By definition, the Riemann curvature tensor $\mathcal{R}$ is an
$\so(4)$-valued 2-form.  However, using the identification
$\Lambda^{2} \cong \so(4)$, we may view it as a self-adjoint
linear map $\mathcal{R} : \Lambda^{2} \rightarrow\Lambda^{2}$
given in coordinates by
\begin{displaymath}
	\mathcal{R}\left(dx^{i}\wedge dx^{j}\right) =
	    \half \sum_{k,l} R_{ijkl}\,dx^{k}\wedge dx^{l}.
\end{displaymath}
Relative to the familar decomposition
$\Lambda^{2} = \Lambda^{2}_{+} \oplus \Lambda^2_{-}$ of the space
of two-forms into its self-dual and anti-self-dual subspaces, the
Riemann curvature can be written in the block matrix form
\begin{displaymath}
	\mathcal{R} = \left( \begin{array}{c|c}
	                      \mathcal{W}^{+} - \frac{1}{12}R & R_{0} \\ \hline
	                      R_{0}^{\ast} & \mathcal{W}^{-} - \frac{1}{12}R
	                  \end{array}
	           \right).
\end{displaymath}
Here $R$ denotes the scalar curvature multiplied by the identity matrix,
while $R_{0} : \Lambda^{2}_{-}\rightarrow \Lambda^{2}_{+}$ is the
trace-free Ricci curvature tensor,
$R_{0}^{\ast} : \Lambda^{2}_{+}\rightarrow \Lambda^{2}_{-}$ is its
adjoint, and $\mathcal{W} = \mathcal{W}^{+} + \mathcal{W}^{-}$ is the conformally invariant Weyl tensor.  A standard reference for this material is 
\cite{AHS}.

We now consider the Riemann curvature $\mathcal{R}'$ of the metric $g'$ 
discussed above.  Since $g'$ is by definition conformally flat, we see
that the Weyl tensor $\mathcal{W}'$ vanishes.  We also recall from
Lemma~\ref{lemma-scalar-curvature} that if our
super-potential $\r$ is harmonic, then the scalar curvature $R'$
vanishes as well.  All that remains is the trace-free Ricci tensor $R_{0}'$,
and so the Riemann curvature is simply
\begin{displaymath}
	\mathcal{R}' = \left( \begin{array}{c|c}
	                      0 & R_{0}'\\ \hline
	                      R_{0}'^{\ast} & 0
	                   \end{array}
	            \right).
\end{displaymath}
Note that the splitting $\spin(4) \cong \sp(1) \oplus \sp(1)$ which we
used to construct the connections $A^{+}$ and $A^{-}$ is isomorphic to
the decomposition $\Lambda^{2} = \Lambda^{2}_{+} \oplus \Lambda^{2}_{-}$.
We can therefore read off the curvatures $F_{A^{+}}$ and $F_{A^{-}}$ of
these connections directly from the block form of the Riemann curvature,
giving us
$$
	F_{A^{+}}^{+} = 0   \qquad
    F_{A^{+}}^{-} = R_{0}' \qquad
	F_{A^{-}}^{+} = R_{0}'^{\ast} \qquad
	F_{A^{-}}^{-} = 0.
$$
Hence the connections $A^{+}$ and $A^{-}$ are anti-self-dual and self-dual
respectively as claimed in Theorem \ref{theorem-1}.


\section{Hyperbolic Vortices}

\subsection{Dimensional Reduction}

\label{dimensional-reduction}

In this section, we examine $\SO(3)$-invariant instantons, showing that
the SD and ASD equations for Yang-Mills connections over $S^{4}$ with $\SO(3)$
symmetry are equivalent to the $\U(1)$ vortex equations over the hyperbolic
plane $\mathcal{H}^{2}$.  This is an example of dimensional reduction, whereby
the Yang-Mills or (A)SD equations for a symmetric connection reduce to
differential equations for a connection and Higgs fields (sections of the
Lie algebra bundle) over a lower dimensional space.

Viewing $S^{4}$ as the standard conformal compactification
$\R^{4} \cup \{\infty\}$, we let $\SO(3)$ act via its fundamental
representation on a three-dimensional subspace of $\R^{4}$.
Expressing this using quaternionic notation, we see that an element
$g\in \Sp(1)$ acts on $\H$ according to $g:x\mapsto g x g^{-1}$, fixing
the real part of $x$ and acting by the adjoint representation on its 
imaginary part. Regarding $S^{4}$ as the quaternionic projective space
$P(\H^{2})$ with homogeneous coordinates $(x:y)=(x\alpha:y\alpha)$,
the embedding of $\H$ is simply the map $x\mapsto(1:x)$.
The $\Sp(1)$-action given by
    \[g:(x:y)\mapsto(gx:gy) = (gxg^{-1}:gyg^{-1})\]
then provides an extension of the above action on $\R^{4}$ to all of $S^{4}$.

In order to discuss $\SO(3)$-invariant connections, we must lift this 
action on $S^{4}$ to an action on the Lie algebra bundle with fibres
$\so(3)$.  There are two possible lifts: either $\SO(3)$ acts trivially
on each fibre or it acts via the adjoint representation.  For our
purposes, we will consider this second, more interesting, action.
Note that for any $g\in\SO(3)$, the adjoint action leaves fixed an
$\R = \u(1)$ subalgebra.


Again adopting quaternionic notation, any $x\in\H$ can be written
in the form $x = t + r Q$, with $t,r$ real, $r \geq 0$, and Q pure 
imaginary with $Q^{2} = -1$.
Note that $t,r$ coordinatize the upper half-plane, which we
will later regard as hyperbolic space $\mathcal{H}^{2}$.
If $A$ is an $\Sp(1)$-invariant connection, then its connection
one-form satisfies $g A(t,r,Q) g^{-1} = A(t,r,gQg^{-1})$.
The most general connection exhibiting this symmetry is of the form
\begin{displaymath}
	A = \half\,\bigl(Qa + \p_{1}\,dQ + \p_{2}\,Q\,dQ\bigr), \label{symmetric-connection}
\end{displaymath}
where $a = a_{t}\,dt + a_{r}\,dr$, and the $a_{t}, a_{r}, \p_{1}, \p_{2}$ 
are all real functions of $t,r$.  The curvature $F_{A}$ of this connection
$A$ is then
\begin{align*}
	F_{A}&=\half \left( Q\,da\,+\,
	        \half \left(\p_{1}^{2}+\p_{2}^{2}+2\p_{2}
	               \right)dQ \wedge dQ\,+\right.\\*
	     & \qquad\qquad\left.\rule{0in}{3ex}
	        [ d\p_{1} - a ( \p_{2} + 1 ) ] \wedge dQ\,+\,
	        ( d\p_{2} + a \p_{1} ) \wedge Q\,dQ \right).
\end{align*}
Putting $\p = \p_{1}+ i\,(\p_{2}+1)$ and writing $d_{a}\p = d\p + ia\p$,
the curvature can be written much more simply as
\begin{align*}
	F_{A}=\half \left( Q\,da\,-\,
	        \half\,\left( 1 - |\p|^{2}\right)\,dQ \wedge dQ\,+
	        \Re(d_{a}\p) \wedge dQ \,+\,
	        \Im(d_{a}\p) \wedge Q\,dQ \right).
\end{align*}
Note that multiplication by $Q$ here behaves like multiplication by $i$.


From the above discussion, we see that an $\SO(3)$-invariant
connection $A$ on $S^{4}$ gives rise in a natural way to a $\U(1)$
connection $ia$ and a complex scalar field $\p$ on the upper half-plane.
The next step is to analyze the SD and ASD equations in terms of this
dimensional reduction.  To determine the action of the Hodge star operator,
we consider the 2-form $dx\wedge d\xbar$ which we already know to be
self-dual.
In coordinates $t,r,Q$, we have
\begin{align*}
	dx \wedge d\xbar &= (dt+Q\,dr+r\,dQ) \wedge (dt-Q\,dr-r\,dQ) \\*
	    &= 2Q\,dt\wedge dr\,+\,r^{2}\,dQ\wedge dQ\,+\,
	        2r\,(dt\wedge dQ\,+\,dr\wedge Q\,dQ),
\end{align*}
and so the Hodge star operator acts according to
\begin{align*}
	\ast\,Q\,dt\wedge dr & =  \frac{r^{2}}{2}\,dQ\wedge dQ \\
	\ast\,dt\wedge dQ    & =  dr\wedge Q\,dQ \\
	\ast\,dr\wedge dQ    & =  -dt\wedge Q\,dQ.
\end{align*}
Furthermore, using the hyperbolic metric
\begin{displaymath}
	h = \frac{1}{r^{2}}\,( dt^{2} + dr^{2} )
\end{displaymath}
 on the upper half-plane, the corresponding Hodge star operator
 $\ast_{h}$ satisfies $\ast_{h}dt \wedge dr = r^{2}$,
 $\ast_{h}dt = dr$, and $\ast_{h}dr = -dt$. Combining this with
 the usual Hodge star operator yields
\begin{align*}
	                   \ast\,Q\,da  =  \half\,(\ast_{h} da)\,dQ \wedge dQ \qquad
	\ast\,\Re\,(d_{a}\p) \wedge dQ  =  \ast_{h} \Re\,(d_{a}\p) \wedge Q\,dQ.
\end{align*}
Using complex notation with $z=t+ir$ and noting the identity
$\ast_{h}dz = -i\,dz$, we observe that
$2\,\dbar_{a}\p = d_{a}\p - i \ast_{h}d_{a}\p$.  We therefore
conclude that the $\SO(3)$-symmetric self-duality equation
$F_{A} = \ast F_{A}$ on $S^{4}$ is equivalent to the following
two equations on hyperbolic space $\mathcal{H}^{2}$:
\begin{align}
	   \dbar_{a} \p & =  0 \label{eq:vortex-a} \\
   \ast_{h} i F_{a} & =   1 - |\p|^{2}, \label{eq:vortex-b}
\end{align}
where $F_{a} = i\,da$ is the curvature of the connection $ia$.
These equations are known as the {\em vortex equations}. Similarly, the
$\SO(3)$-symmetric anti-self-dual equation  $F_{A} = -\ast F_{A}$
is equivalent to the {\em anti-vortex equations}:
\begin{align}
    \partial_{a} \p & =  0 \label{eq:anti-vortex-a} \\
   \ast_{h} i F_{a} & =  |\p|^{2} - 1. \label{eq:anti-vortex-b}
\end{align}
The first equation in each pair is simply the condition that $\p$ be
holomorphic (or anti-holomorphic) with respect to the holomorphic
structure compatible with the connection $ia$.  The second equation
then expresses a form of duality between the connection and Higgs field.
These vortex equations are discussed in great detail in \cite{JT}%
\SSfootnote{After adjusting to the slightly different notation of
          \cite{JT}, using $a' = -a$ and $\p' = \bar{\p}$, the
          reader will find that equations (11.5a) and (11.5b) on
          p. 99 of \cite{JT} should be switched.}.
Note that if we consider these vortex and anti-vortex equations over
the plane $\R^{2}$ with the flat metric $h = \half(dx^{2}+dy^{2})$,
then we obtain the Euclidean vortex and anti-vortex equations in their
customary form as given by (1.7) and (1.8) on p. 55 of \cite{JT}. 


We now compute the $L^{2}$ norm of the curvature $F_{A}$ using the 
standard metric $|\xi|^{2} = \xi\bar{\xi} = -\xi^{2}$ on the Lie
algebra $\sp(1)$ of imaginary quaterions. The Yang-Mills action,
or energy, of this $\SO(3)$-invariant connection is thus
\begin{align*}
    \|F_{A}\|^{2} & =  \int_{S^{4}}
                        -\,F_{A}\wedge \ast F_{A} \\
                  & =  \frac{1}{8}\int_{S^{4}}
              \Bigl( da\wedge \ast_{h} da
        \,+\, (1-|\p|^{2}) \wedge \ast_{h}(1-|\p|^{2}) \,+\,
              \\
    & \qquad \qquad  \rule{0in}{2ex}
              2\,\Re d_{a}\p \wedge\ast_{h}\Re d_{a}\p
              \,+\,2\,\Im d_{a}\p \wedge\ast_{h}\Im d_{a}\p
                    \Bigr) 
 	\wedge\,\left(\,dQ\wedge Q\,dQ\,\right).
\end{align*}
Note that the left factor of the integrand is independent of the
variable $Q$.  Since $Q$ parametrizes the unit 2-sphere with volume
form $\half(dQ\wedge Q\,dQ)$, we can integrate out a factor of
$\int_{S^{2}}dQ\wedge Q\,dQ = 8\pi$, leaving an integral over
hyperbolic space.  The action then becomes
\begin{displaymath}
   \|F_{A}\|^{2} = \pi \left( \|F_{a}\|^{2}_{h} \,+\,
                              2\,\|d_{a}\p\|^{2}_{h} \,+\,
                              \|1-|\p|^{2}\|^{2}_{h}
                       \right),
\end{displaymath}
which we recognize as the $\U(1)$ Yang-Mills-Higgs action on hyperbolic
space, at least up to a constant. From this action, we see that a
finite-energy $\SO(3)$-invariant connection on $S^{4}$ corresponds to a
pair $(ia,\p)$ over $\mathcal{H}^2$ satisfying the boundary conditions
\begin{displaymath}
	d_{a}\p(x) \rightarrow 0, \qquad |\p(x)| \rightarrow 1,
\end{displaymath}
as $|x| \rightarrow \infty$.


Next we examine the relationship between the Chern classes of an
$\SO(3)$-invariant connection $A$ on $S^{4}$ and those of the corresponding
connection $ia$ over $\mathcal{H}^2$.  Computing $c_{2}(A)$, we first
note that the negative definite form $\xi\mapsto\mbox{Tr}(\xi)^2$
on the Lie algebra $\su(2)$ corresponds to $\xi\mapsto 2\xi^{2}$
on $\sp(1)$. In this quaternionic notation we therefore have
\begin{align*}
    c_{2}(A) & =  -\frac{1}{4\pi^{2}}
                   \int_{S^{4}}F_{A}\wedge F_{A} \\
             & =  -\frac{1}{16\pi^{2}}
                   \int_{S^{4}} \left[\,
                       ( 1-|\p|^{2} )\,da \,-\,
                       2\,\Re d_{a}\p\wedge\Im d_{a}\p
                                \,\right] 
                \wedge\,
                   \left(\,dQ\wedge Q\,dQ\,\right) \\
             & =  \frac{i}{2\pi} \int_{\mathcal{H}^2} F_{a}
             =  c_{1}(a).
\end{align*}
Here we again integrate out the $S^{2}$ factor $dQ\wedge Q\,dQ$, and
on the last line we apply Stokes' theorem with the integrand
\begin{align*}
	d\,(i\bar{\p}\,d\p) & =  i\,d\bar{\p} \wedge d\p \,-\,
	                        \p\bar{\p}\,da \,-\,
	                        (\p\,d\bar{\p} + \bar{\p}d\p)\wedge a \\*
	                  & =  -|\p|^2\,da - 2\,\Re d_{a}\p\wedge\Im d_{a}\p,
\end{align*}
assuming that $i\bar{\p}\,d\p$ vanishes at infinity.


\subsection{Another Harmonic Function ansatz}
\label{vortex-ansatz}

We now return to the harmonic function ansatz that we discussed in 
Section~\ref{instanton-ansatz}.  If we begin with a $\SO(3)$-invariant 
harmonic super-potential, then the resulting SD or ASD connection will
also exhibit $\SO(3)$ symmetry, and from the previous section
we know that such a connection is equivalent to a hyperbolic vortex
or anti-vortex.
In this section, we take advantage of this dimensional reduction
to provide a similar harmonic function ansatz constructing solutions to
the vortex equation over hyperbolic space.


Our first step is to examine the relationship between the Laplacian
on hyperbolic space $\mathcal{H}^{2}$ and the $\SO(3)$-symmetric Laplacian
on $\R^{4}$.

\begin{lemma}
    An $\SO(3)$-invariant function $\r$ on $\R^{4}$ is harmonic if and
    only if it can be written as $\r = r^{-1}\phi$, where $\phi$ is
    a harmonic function on $\mathcal{H}^2$.
\end{lemma}

\begin{proof}
Using the quaternionic notation $x = t + rQ$ introduced in
\S\ref{dimensional-reduction}, if $\r = \r(r,t)$ is an $\SO(3)$-invariant
function on $\R^{4}$, then its Laplacian is
\begin{equation} \label{eq:laplacian}
	\Delta \r = -\left( \delsq{t} + \delsq{r} +
	                    \frac{2}{r}\del{r}
	             \right) \r.       
\end{equation}
To cancel the unwanted linear term, we put $\r = r^{-1}\phi$,
where $\phi$ is a function on $\mathcal{H}^{2}$. The Laplacian $\Delta$
then becomes
\begin{displaymath}
	\Delta r^{-1}\phi = -\frac{1}{r} \left( \delsq{t} + \delsq{r}
	                                 \right) \phi
	                  = r^{-3}\Delta_{h} \phi,
\end{displaymath}
where the Laplacian $\Delta_{h}$ on $\mathcal{H}^{2}$ is
\begin{equation} \label{eq:hyperbolic-laplacian}
    \Delta_{h} = -r^{2}\left(\delsq{t} + \delsq{r}\right)
\end{equation}
These two Laplacians are therefore related by%
\SSfootnote{In general, the conformal Laplacian is $L_{g}=\Delta+kR$,
where $R$ is the scalar curvature and $k$ is a constant depending on
the dimension. Taking the metric $g' = e^{2f}g$, it we have
\begin{displaymath}
    L_{g'} = e^{-(d+2)f/2}L_{g}\,e^{(d-2)f/2}.
\end{displaymath}
Here the Euclidean metric on $R^{4}$ is $g=dt^{2}+dr^{2}+r^{2}dS^{2}$,
where $dS^{2}$ is the metric on $S^{2}$. Taking the conformally equivalent
metric $g' = r^{-2}(dt^{2}+dr^{2}) + dS^{2}$ on $\mathcal{H}^2 \times S^{2}$,
Lemma \ref{lemma-scalar-curvature} tells us that $R' = 0$ since $r^{-1}$
is harmonic. It follows that $\Delta'=r^{3}\Delta r^{-1}$.}
$\Delta$ = $r^{-3}\Delta_{h}r$ and thus $\r$ is harmonic on
$\R^{4}$ if and only if $\phi$ is harmonic on $\mathcal{H}^{2}$.
\end{proof}


We recall from Theorem \ref{theorem-1} that the connection given by
equation~(\ref{eq:self-dual}),
\begin{displaymath}
    A = - \Im\left(\del{x}\log \r\,dx \right),
\end{displaymath}
is self-dual if and only if $\r$ is harmonic. In our current notation,
the quaternionic differential and partial derivative in this expression are
\begin{align*}
	dx       =  dt + Q\,dr + r\,dQ \qquad
	\del{x}  =  \half \left( \del{t} - Q\,\del{r} - \cdots
	                        \right),
\end{align*}
where we have left out the portions of the partial derivative in the $Q$
directions as these vanish when applied to $\SO(3)$-invariant functions.
Taking $\r = r^{-1}\phi$, we see that $\log \r = \log \phi - \log r$.  
Expanding equation~(\ref{eq:self-dual}) using these expressions,
our $\SO(3)$-invariant self-dual connection becomes
\begin{align*}
    A & =  \half\,\left(\,Q \left[\,
                  \left( \del{r}\log\phi - \frac{1}{r}\right) dt \,-\,
                  \del{t}\log\phi\,dr\,
                           \right]\right. \\*
      &   \qquad \left. \mbox{\qquad}
                  -\,r\,\del{t}\log\phi\,dQ\,+\,
                  \left( r\,\del{r}\log\phi - 1 \right) Q\,dQ\,
            \right).
\end{align*}
As we did in in \S\ref{dimensional-reduction}, we can extract from
this connection the $\U(1)$ connection
\begin{equation} \label{eq:vortex-da}
	d_a = d\,+\,i \left[\,
	          \left( \del{r}\log\phi - \frac{1}{r}\right) dt -
              \del{t}\log\phi\,dr\,
                  \right]
\end{equation}
with curvature
\begin{align*}
	F_{a}  =  -i\left[
	         \left(\delsq{t} + \delsq{r}\right) \log\phi +
	         \frac{1}{r^{2}}
	              \right] dt\wedge dr
	       =  -i \left( 1 - \Delta_{h}\log\phi
	                 \right)\,r^{-2}dt\wedge dr,
\end{align*}
and the complex Higgs field
\begin{equation} \label{eq:vortex-phi}
    \p = r \left( - \del{t}\log\phi + i\,\del{r}\log\phi \right)    
\end{equation}
with norm
\begin{displaymath}
    |\p|^{2} = r^{2}\left[
	          \left(\del{t}\log\phi\right)^{2} +
	          \left(\del{r}\log\phi\right)^{2}
	                \right] \\*
	         = \left| \nabla\log\phi \right|_{h}^{2}.
\end{displaymath}
Writing the pair $(a,\phi)$ using complex notation with $z = t + ir$,
we obtain the hyperbolic space analogue of Theorem \ref{theorem-1}.

\begin{theorem} \label{theorem-vortex-ansatz}
	Given a positive real-valued super-potential $\phi$ on the hyperbolic
	upper half-plane $\mathcal{H}^{2}$, the connection and Higgs field pair
	$(a,\p)$ defined by
	\begin{align}
	    \dbar_{a}  =  \dbar\,+\,\dbar\log\phi\,+\,\frac{d\zbar}{z-\zbar}
	        \label{eq:vortex-dbar} \qquad
	    \p         =  i\,(z-\zbar)\,\del{z}\log\phi,
	\end{align}
	satisfies the vortex equations (\ref{eq:vortex-a}) and 
	(\ref{eq:vortex-b}) and the pair $(a',\p')$ defined by
	\begin{align} \label{eq:anti-vortex-dbar}
	    \partial_{a}' & =  \partial\,+\,\partial\log\phi\,-\,\frac{dz}{z-\zbar} \\
	    \p'           & =  -i\,(z-\zbar)\,\del{\zbar}\log\phi,
	                        \label{eq:anti-vortex-field}
	\end{align}
	satisfies the anti-vortex equations (\ref{eq:anti-vortex-a})
	and (\ref{eq:anti-vortex-b}) if and only if the super-potential
	$\phi$ is harmonic.
\end{theorem}

\begin{proof}
Recalling that $\ast_{h} 1 = r^{-2}dt\wedge dr$ with our
hyperbolic metric, we see that the second of the vortex equations
$iF_{a} = \ast_{h}\left( 1 - |\p|^{2} \right)$ reduces to
\begin{displaymath}
    \frac{1}{\phi}\,\Delta_{h}\phi
      = \Delta_{h}\log\phi - \left| \nabla\log\phi \right|_{h}^{2}
	  = 0.
\end{displaymath}
Using the complex form (\ref{eq:vortex-dbar}),
it is easy to verify that the holomorphicity condition $\dbar_{a}\p = 0$
likewise reduces to $\Delta_{h} \phi = 0$. Hence the pair $(a, \p)$ satisfies
the vortex equations if and only if the super-potential $\phi$ is harmonic.
Similarly, the pair $(a',\phi')$ is derived by dimensional reduction from
the ASD connection (\ref{eq:anti-self-dual}), and so it satisfies the
anti-vortex equations if and only if $\phi$ is harmonic.
\end{proof}


Computing the Chern class $c_{1}$ for a pair $(a,\p)$ satisfying
the second of the vortex equations (\ref{eq:vortex-b}), we have
\begin{displaymath}
	c_{1}(a) = \frac{i}{2\pi}\int_{\mathcal{H}^{2}} F_{a}
	         = \frac{1}{2\pi}\int_{\mathcal{H}^{2}}
	               \ast_{h}\left(1-|\p|^{2}\right).
\end{displaymath}
For the vortex over the upper half-plane $\mathcal{H}^{2}$ constructed
in (\ref{eq:vortex-da}) and  (\ref{eq:vortex-phi}) using a
harmonic super-potential $\phi$, this Chern class takes the form
\begin{align*}
	c_{1}(a)=\frac{1}{2\pi}\int_{\R^{2}_{+}}\left(
	               1 - \Delta_{h}\log\phi
	           \right)\,r^{-2} dt\wedge dr
	        =\frac{1}{2\pi}\int_{\R^{2}_{+}}\left(
	               \frac{1}{r^{2}} - 4\left|\del{z}\log\phi\right|^{2}
	           \right)\,dt\wedge dr.
\end{align*}
Likewise, if we use the above ansatz to construct the anti-vortex
corresponding to a harmonic super-potential, then the Chern class
switches sign.


As in \S\ref{quaternionic-notation}, we note that the curvature
$F_{a}$ of a vortex is gauge invariant.  Therefore, if two harmonic
super-potentials $\phi_{1}$ and $\phi_{2}$ over hyperbolic space
yield gauge equivalent vortices, then they must satisfy the equations
\begin{displaymath}
	\left|\del{z}\log\phi_{1}\right| = \left|\del{z}\log\phi_{2}\right|
\end{displaymath}
and $\Delta_{h}\log\phi_{1} = \Delta_{h}\log\phi_{2}$.


Note that it is significantly simpler to calculate $c_{1}$ directly
from the vortex construction on $\mathcal{H}^{2}$ than it is by invoking 
dimensional reduction and computing the equivalent Chern class $c_{2}$
for the corresponding $\SO(3)$-invariant instanton over $S^{4}$. Indeed,
by comparing the above expression for $c_{1}(a)$ with the expression
(\ref{eq:c2}) for $c_{2}(A)$, we obtain a circuitous proof
of the identity
\begin{displaymath}
	\int_{\R^{2}_{+}} \left( 1 - \Delta_{h}\log\phi
	                  \right)\,r^{-2}dt \wedge dr =            
	- \int_{\R^{2}_{+}} \left(
	    \frac{1}{2}\,r^{2} \Delta\Delta\log \frac{\phi}{r}
	                    \right) dt \wedge dr
\end{displaymath}
for a harmonic function $\phi$ defined on the upper half-plane, where
$\Delta_{h}$ is the Laplacian on $\mathcal{H}^{2}$ given by 
(\ref{eq:hyperbolic-laplacian}) and $\Delta$ is the Laplacian on $\R^{4}$
given by (\ref{eq:laplacian}).


\subsection{Conformal Transformations Revisited}

Instead of relying on dimensional reduction to derive the vortex ansatz
of Theorem~\ref{theorem-vortex-ansatz}, we present here an interpretation
of this construction that is entirely intrinsic to hyperbolic space. As
we did in \S\ref{conformal-instanton}, we can treat the super-potential
as a conformal transformation and then compute the Levi-Civita
connection of the resulting metric. Since we are working on
two-dimensional hyperbolic space, we can take advantage of complex
notation to simplify our task.

Let $\dbar$ be the standard holomorphic stucture on the complex upper
half-plane. Choosing a holomorphic tangent frame (i.e., a single
holomorphic section) $e$, consider the Hermitian metric $g$ specified
by $(e,e)_{g} = \r^{2}$, where $\r$ is a smooth nonzero real-valued
function. With respect to our holomorphic frame $e$, the unique
connection compatible with both the holomorphic structure $\dbar$ and
the metric $g$ is specified by the $(1,0)$-form
\begin{displaymath}
    a = \r^{-2}(\partial \r^{2}) = 2\,\partial\log\r.
\end{displaymath}
To express this in the form of a unitary connection (in this case given
by a purely imaginary complex 1-form), we must switch to a tangent frame
that is orthonormal with respect to the metric $g$. In terms of the unitary
frame $e' = \r^{-1}e$, the connection $a$ then becomes
\begin{displaymath}
    a' = \partial\log\r - \dbar\log\r
       = 2i\,\Im\partial\log\r =  -2i\,\Im\dbar\log\r
\end{displaymath}
and the new holomorphic structure is $\dbar'=\dbar-\dbar\log\r$,
which we observe is compatible with the connection $a'$. In either
frame, the curvature of this connection is given by
\begin{displaymath}
	F_{a} = 2\,\dbar\partial\log\r = -i\,\Delta\log\r\,d\mu,
\end{displaymath}
where the volume form $d\mu$ and the Laplacian $\Delta$ are both
taken here with respect to the Euclidean metric on $\R^{2}$.

When working with $\mathcal{H}^2$, the hyperbolic metric $h$ on the upper
half-plane corresponds to the function $\r = r^{-2}$. Taking a conformal 
transformation, we consider the metric $h'$ specified by a function of
the form $\r = \phi^{2}/r^{2}$ with $\phi$ harmonic. The resulting 
unitary connection $a$ then splits into the $(0,1)$ and $(1,0)$ components
\begin{equation*}
	\dbar_{a}     =  \dbar\,-\,\dbar\log\phi\,-\,\frac{d \zbar}{2ir} \qquad
	\partial_{a}  =  \partial \,+\,\partial\log\phi\,-\,\frac{dz}{2ir},
\end{equation*}
noting that $r = (z - \zbar)/2i$. The curvature of this connection is then
\begin{displaymath}
	i F_{a} = \left( \Delta\log\phi - r^{-2} \right) d\mu
	        = \left( \Delta_{h}\log\phi - 1 \right) d\mu_{h}
\end{displaymath}
where $d\mu_{h} = r^{-2}d\mu$ is the volume form and $\Delta_{h} = 
r^{2}\Delta$ is the Laplacian for the hyperbolic metric $h$.
It is then easy to show that the complex Higgs field $\p$ defined by
\begin{displaymath}
	\p = 2r\,\del{\zbar}\log\phi = 2\,\frac{r}{\phi}\,\del{\zbar}\phi
\end{displaymath}
satisfies the anti-vortex equations $\partial_{a}\p = 0$ and
$\ast_{h} i F_{a} = |\p|^{2} - 1$ if the super-potential $\phi$ is
harmonic. We observe that this pair $(a,\p)$ agrees with the anti-vortex
(\ref{eq:anti-vortex-dbar}) and (\ref{eq:anti-vortex-field}) constructed
by dimensional reduction of an anti-self-dual connection over the 4-sphere.
Similarly, we can construct the vortex given by (\ref{eq:vortex-dbar})
by reversing orientation, thereby exchanging
the holomorphic and anti-holomorphic structures $\dbar$ and $\partial$.


\subsection{The Symmetric 't Hooft Construction}

\label{vortex-tHooft}

In \S\ref{tHooft}, as an illustration of the harmonic function ansatz,
we constructed the 't Hooft instantons.  These are the instantons formed
by taking the superposition of multiple copies of the basic instanton
with varying scales and distinct centers.  For our super-potential, we used
a sum of the Green's functions of the Laplacian, centered at the given 
points and weighted according to the corresponding scales.  If we impose 
$\SO(3)$ symmetry on this class of instantons, we see that all of the 
centers must lie on a single real line.  In fact, as we will demonstrate
in the following section, {\em all\/} $\SO(3)$-invariant instantons can be
constructed in this manner---as the superposition of basic instantons on a
line.  In this section, we examine the hyperbolic vortices associated
to these symmetric 't Hooft instantons by dimensional reduction.


We begin with the basic instanton with unit scale centered at the origin,
which we recall is given by the $\R^{4}$ super-potential $\r = 1 + |x|^{-2}$.
The corresponding super-potential for hyperbolic space $\mathcal{H}^{2}$ is then
\begin{equation}\begin{split}
	\phi  =  r\r = r + \frac{r}{r^{2}+t^{2}}  
	      =  \Im \left( z - \frac{1}{z} \right) 
	           = \frac{ (z-\zbar)\,(1 + z\zbar) }{2i\,z\zbar }.
	           \label{eq:simple-potential}
\end{split}\end{equation}
Taking the complex partial derivatives of its logarithm, we obtain
\begin{align*}
    \del{z}\log\phi  =  + \frac{1}{z-\zbar}\,+\,\frac{\zbar}{1+z\zbar}
                          \,-\,\frac{1}{z} \qquad
    \del{\zbar}\log\phi  =  - \frac{1}{z-\zbar}\,+\,\frac{z}{1+z\zbar}
                              \,-\,\frac{1}{\zbar}.
\end{align*}
Inserting these expressions into the formula (\ref{eq:vortex-dbar}) gives us the connection and Higgs field
\begin{equation*}
	\dbar_{a} =  \dbar\,-\,\frac{d\zbar}{\zbar\,(1+z\zbar)} \qquad
	\p        =  i\,\frac{\zbar\,(1+z^{2})}{z\,(1+z\zbar)},
\end{equation*}
satisfying the vortex equations (\ref{eq:vortex-a}) and (\ref{eq:vortex-b}).  Note that
as $z$ approaches the real axis, the Higgs field obeys the
boundary condition $|\p|\rightarrow 1$.


One of the primary results concerning solutions to the vortex equations 
is that they are uniquely specified up to gauge equivalence by the zeros 
of the Higgs field (see \cite[Chapter III]{JT}).
In the example above, we see that $\p$ vanishes at the point $z = i$.
If we alter the scale of our basic instanton and translate it along the 
real axis, the $\mathcal{H}^{2}$ super-potential becomes
\begin{displaymath}
	\phi = \Im \left( z - \frac{\lambda}{z-a} \right)
\end{displaymath}
with $\lambda > 0$ and $a$ real. The corresponding vortex is then
\begin{equation*}
	\dbar_{a}  =  \dbar\,-\,\frac{\lambda\,d\zbar}
	                               {(\zbar-\bar{a})\,(1+|z-a|^{2})} \qquad
	\p         =  i\,\frac{(\zbar-\bar{a})\,(1+(z-a)^{2})}
	                        {(z-a)\,(1+|z-a|^{2})},
\end{equation*}
and we see that $\p$ vanishes at the point $z = a + i\sqrt{\lambda}$.
Most generally, given a set of $k$ complex points $\{z_{i}\}$ in the
upper half-plane, the super-potential
\begin{displaymath}
	\phi = \Im \left( z - \sum_{i = 1}^{k}
	                       \frac{(\Im z_{i})^2}{z - \Re z_{i}}
	           \right)
\end{displaymath}
generates the unique hyperbolic vortex with Higgs field vanishing at
the points $\{z_{i}\}$.  Hence the centers of the instantons correspond 
to the real parts of the complex zeros, while the scales correspond to
their imaginary parts.


\subsection{The Equivariant ADHM Construction}

In this section we shall use an $\SO(3)$ equivariant version of the ADHM 
construction \cite{ADHM} in order to provide an alternative construction for the
symmetric 't Hooft instantons discussed in the previous section.
In addition, since the ADHM construction actually generates {\em all}
possible anti-self-dual connections on bundles over $S^{4}$, we
will then be able to show that every symmetric instanton must be
gauge equivalent to one constructed using the 't Hooft ansatz. By 
dimensional reduction, this gives us a complete classification of 
hyperbolic vortices, proving that such vortices are uniquely determined 
up to a gauge transformation by the zeros of their Higgs fields. This is
to be contrasted with Euclidean vortices, in which case the classification
theorem may be proved using approximation techniques (see \cite[Chapter 
III]{JT}), but no explicit construction for the vortex solutions is known.

\newcommand{\WH}{W\!\otimes\!\H_{1}}

Here we use the quaternionic version of the construction as discussed in
\cite{A}. When dealing with quaternionic vector spaces and linear maps,
we use the convention that scalar multiplication acts on the {\em right}.
Recall from \S\ref{dimensional-reduction} that under our $\Sp(1)$-action,
we may view $S^{4}$ as the quaternionic projective space $P(\H_{1}^{2})$,
where $\H_{1}$ is the fundamental representation with $\Sp(1)$ acting by 
{\em left} quaternion multiplication.

To construct an $\Sp(1)$-invariant ASD connection on the bundle
$E\rightarrow S^{4}$ with Chern class $c_{2}(E) = -k$, we introduce
the $k+1$ dimensional {\em quaternionic\/} $\Sp(1)$ representation $V$
given by
\begin{equation} \label{eq:V}
	V = \mbox{Ker}\,\mathcal{D}_{A}^{\ast} :
	        \Gamma (S^{4}, E\otimes S^{-}\otimes S^{-})
	            \rightarrow
	        \Gamma (S^{4}, E\otimes S^{+}\otimes S^{-})
\end{equation}
and the $k$ dimensional {\em real\/} $\Sp(1)$ representation $W$ given by
\begin{equation} \label{eq:W}
    W = \left( \mbox{Ker}\,\mathcal{D}_{A}^{\ast} :
	        \Gamma (S^{4}, E\otimes S^{-})
	            \rightarrow
	        \Gamma (S^{4}, E\otimes S^{+})
        \right)_{\R}^{\ast},
\end{equation}
where $A$ is an arbitrary connection on $E$ (the spaces $V,W$ are 
independent of the connection), $S^{\pm}$ are the two quaternionic
half-spin bundles, and $\mathcal{D}_{A}^{\ast}$ is the adjoint of the Dirac
operator with coefficients in $E\otimes S^{-}$ and $E$ respectively.
Using these spaces $V,W$ the ADHM data consists of the the three maps:
\begin{itemize} \samepage
    \item an arbitrary $\Sp(1)$ equivariant inclusion
          $\WH \hookrightarrow V$
    \item an $\Sp(1)$ equivariant $\H$-linear map
          $B:\WH\rightarrow\WH$ satisfying $B^{\ast} = \bar{B}$
          (i.e., $B$ is represented by a symmetric matrix)
    \item an $\Sp(1)$ equivariant $\H$-linear map
          $\Lambda : \WH \rightarrow V\,/\,\WH$.
\end{itemize}
If we fix the inclusion $\WH \hookrightarrow V$, then we say that
two sets of ADHM data $(B,\Lambda)$ and $(B',\Lambda')$ are
{\em equivalent\/} if
\begin{displaymath}
	B' = U B U^{-1}, \qquad  \Lambda' = v \Lambda U^{-1}
\end{displaymath}
for suitable $U\in O(W)$ and $v\in \Sp(V\,/\,\WH)$. From the ADHM data,
we construct an $\Sp(1)$ equivariant family of $\H$-linear maps
$v(x) : \WH \rightarrow V$ parametrized by $x\in\H$, given by
\begin{displaymath}
	v(x) = \left( \begin{array}{c}
	           \Lambda \\
	           B - xI
	       \end{array} \right)
\end{displaymath}
relative to the decomposition $V = (V\,/\,\WH) \oplus (\WH)$.

\begin{theorem}[ADHM] \label{ADHM}
    There is a one-to-one correspondence between equivalence classes of
    ADHM data $(B,\Lambda)$ satisfying the two conditions
    \begin{description}
        \item[non-degeneracy] $v(x)$ is injective for all $x\in\H$
        
        \item[ADHM condition] $\Lambda^{\ast}\Lambda + B^{\ast}B :
                              \WH \rightarrow \WH$ is real,
    \end{description}
    and gauge equivalence classes of $\Sp(1)$-invariant ASD connections
    on $E$.
\end{theorem}    

To construct the connection associated to a set of ADHM data, we
first observe that the non-degeneracy condition implies that
$f : x \mapsto \mbox{Coker}\,v(x)\subset V$ is a smooth map from
$S^{4}$ to the quaternionic projective space $P(V)$ (we map the point
at $\infty$ to the line $V\,/\,\WH$). We then define the corresponding
vector bundle $E$ and connection $A$ to be the pullback of the canonical
quaternionic line bundle over $P(V)$ with its standard connection
(induced by orthogonal projection from the trivial flat connection on V).
The fact that $A$ is ASD follows from the ADHM condition. For a complete
proof of the non-equivariant version of this theorem, see 
\cite[\S3.3]{DK} or \cite{A}. The proof of equivariant version then
proceeds with minimal modification.

We now give an even more precise description of $\Sp(1)$-invariant 
instantons, starting by examining the characters of the representations
$V$ and $W$.

\begin{lemma}
    The ADHM representations $V,W$ described in (\ref{eq:V}) and (\ref{eq:W})
    corresponding to the bundle $E \rightarrow S^{4}$ with $c_{2}(E) = -k$
    are given by
    \begin{displaymath}
        V = \H_{1}^{k+1}, \qquad W = \R_{0}^{k},
    \end{displaymath}
    where $\H_{1}$ is the fundamental representation of $\Sp(1)$ acting
    by left multiplication and $\R_{0}$ is the trivial real representation.
\end{lemma}

\begin{proof}
If $E$ is an $\Sp(1)$ equivariant vector bundle, then it is
clearly also equivariant with respect to any one-parameter subgroup $S^{1}$ of
$\Sp(1)$. We may therefore apply the results of \cite{Br} and \cite{BA}.
The fixed point set for this $S^{1}$-action on $S^{4}$ is the sphere
$S^{2}$, over which the bundle $E$ splits as $E|_{S^{2}}=L\oplus L^{\ast}$.
Here $L$ is a complex line bundle with an $S^{1}$-action and $L^{\ast}$ 
is its dual. As in \cite{Br}, such $S^{1}$ equivariant bundles $E$ are 
characterized by a pair of constants $(m,k)$, where $m$ is the weight of 
the $S^{1}$-action on $L$ and $k = c_{1}(L^{\ast})$. In our case this 
$S^{1}$-action is derived from the fundamental representation of  $\Sp(1)$,
and so its weight is $m = 1/2$. In addition, noting that $2mk = -c_{2}(E)$,
we see that the two definitions of $k$ agree. Using the equivariant index
calculations of \cite{Br}, Braam and Austin compute the representations $V, W$ in
equations (3.4) and (3.5) of \cite{BA}. For $m = 1/2$ we have
$V=(\C_{1/2}\oplus\C_{-1/2})^{k+1}$ and $W=\R_{0}^{k}$ as
representations of $S^{1}$. The corresponding $\Sp(1)$ representations
are then $V = \H_{1}^{k+1}$ and $W = \R_{0}^{k}$.
\end{proof}

From this lemma, we see that the ADHM data $(B,\Lambda)$ is a pair of
$\Sp(1)$ equivariant maps $B : \H_{1}^{k}\rightarrow\H_{1}^{k}$
and $\Lambda : \H_{1}\rightarrow\H_{1}$. We then have
$gBg^{-1} = B$ and $g\Lambda g^{-1} = \Lambda$ for all $g\in\Sp(1)$,
where $g$ acts by quaternion multiplication. It follows that both
$B$ and $\Lambda$ must be {\em real\/} transformations. Recalling that
$B$ is symmetric, we see that $B$ is diagonalizable with real
eigenvalues. Choosing a suitable basis, we can therefore write the map
$v(x) : \H_{1}^{k} \rightarrow \H_{1}^{k+1}$ as a matrix of the form
\begin{displaymath}
	v(x) = \left( \begin{array}{ccc}
	                  \lambda_{1} & \cdots & \lambda_{k} \\
	                  b_{1} - x   &    0    &   0        \\
	                    0         & \ddots  &   0        \\
	                    0         &   0     & b_{k} - x
	              \end{array}
	       \right)
\end{displaymath}
with real eigenvalues $\{b_{1},\ldots,b_{k}\}$ and real
scales $\{\lambda_{1},\ldots,\lambda_{k}\}$ with $\lambda_{i} \geq 0$.
The non-degeneracy condition of Theorem \ref{ADHM} implies that
the $b_{i}$ are distinct and the $\lambda_{i}$ are nonzero, while the
ADHM condition is automatically satisfied since $B$ and $\Lambda$ are
real. Computing $\mbox{Coker}\,v(x) = \mbox{Ker}\,v(x)^{\ast}$, we obtain

\begin{theorem}
    Given $k$ distinct real centers $\{b_{1},\ldots,b_{k}\}$ and
    $k$ positive real scales $\{\lambda_{1},\ldots,\lambda_{k}\}$, let
    $f : S^{4} \rightarrow P(\H_{1}^{k+1})$ be the map given by
    \begin{displaymath}
		f\,:\,x\,\mapsto\,\left( 1\,:\,\frac{\lambda_{1}}{x-b_{1}}\,:\,
		                    \cdots\,:\,\frac{\lambda_{k}}{x-b_{k}} \right).
	\end{displaymath}
    The $\Sp(1)$-invariant connection obtained by taking the pullback
    of the standard connection on the canonical bundle over $P(\H_{1}^{k+1})$
    is then ASD and has Chern class $c_{2} = -k$.  Furthermore, 
	every $\Sp(1)$-invariant connection on a bundle $E\rightarrow S^{4}$
	with $c_{2}(E) = -k$ is gauge equivalent to one of this form.
\end{theorem}

We note that the instantons constructed by the above theorem are the
superposition of $k$ basic instantons with distinct centers along the
real axis, which we recognize as the 't Hooft instantons from the previous
section in a different guise. By dimensional reduction, we therefore obtain
a constructive proof of the classification theorem for hyperbolic vortices.

\begin{corollary}
    Given a set of $k$ distinct points $\{z_{1},\ldots,z_{k}\}$ in the
    complex upper half-plane, there exists a finite action solution
    $(a,\p)$ to the hyperbolic vortex equations, unique up to gauge
    transformation, such that $\{z_{1},\ldots,z_{k}\}$ is the zero set
    of the Higgs field $\p$. The Chern class of such a vortex is then
    $c_{1}(a) = k$, obtained by counting the zeros of the Higgs field.
\end{corollary}


\subsection{Symmetric Gauge Transformations}

\label{gauge-transformation}

Now that we understand the relationship between symmetric instantons and 
hyperbolic vortices, we examine how the notion of gauge equivalence 
behaves under this dimensional reduction.  Starting with a symmetric
$\SU(2)$ gauge transformation on $S^{4}$, we compute the resulting
$\U(1)$ gauge transformation on hyperbolic space $\mathcal{H}^{2}$.  Then,
continuing to work in the simpler $\mathcal{H}^{2}$ picture, we examine 
the conditions under which two hyperbolic vortices given by the
harmonic function ansatz of \S\ref{vortex-ansatz} are gauge
equivalent.


Using the quaternionic notation $x = t + rQ$, the most general 
$\Sp(1)$-invariant gauge transformation on $S^{4}$ has the form
\begin{displaymath}
	g = e^{Q\,\x(t,r)} = \cos\x(t,r) + Q\sin\x(t,r),
\end{displaymath}
where $\x$ is a real-valued function on $\mathcal{H}^{2}$. Computing
its differential, we have
\begin{equation*}
	dg\,g^{-1}  =  Q\,d\x\,+\,\sin\x\:dQ\:e^{-Q\x} 
	            =  Q\,d\x\,-\,\half\,Q\,(e^{2Q\x}-1)\,dQ.
\end{equation*}
Applying this gauge transformation to the general $\SO(3)$-invariant
connection (\ref{symmetric-connection}), we obtain
\begin{align*}
	g(A)&= \half \, e^{Q\x}\,\left[\,
	              Qa\,+\,(\p_{1}+Q\p_{2})\,dQ
	              \,\right]\,e^{-Q\x}
	        \,-\,dg\,g^{-1} \\
	    &=\half\left[\,
	       Q\,(a - 2\,d\x) \,+\,
	       \left( e^{2Q\x}\,[\,\p_{1}+Q\,(\p_{2}+1)\,] - Q \right)dQ
	              \,\right],
\end{align*}
noting that $g\,Q = Q\,g$ while $g\,dQ = dQ\,g^{-1}$.
The associated connection $ia$ and Higgs field $\p = \p_{1} + 
i\,(\p_{2}+1)$ over $\mathcal{H}^{2}$ then transform according to
\begin{displaymath}
	g(a)   =  a - 2\,d\x, \qquad
	g(\p)  =  e^{2i\x}\p.
\end{displaymath}
Hence, the corresponding gauge transformation over $\mathcal{H}^{2}$
is simply $g = e^{2i\x}$.


Suppose that we have two gauge equivalent hyperbolic vortices
$(a_{+},\p_{+})$ and $(a_{-},\p_{-})$ constructed by the
harmonic function ansatz of \S\ref{vortex-ansatz}, using
equation (\ref{eq:vortex-dbar})
with the super-potentials $\phi_{+}$ and $\phi_{-}$ respectively.
If these two vortices satisfy $a_{-} = g(a_{+})$ and
$\p_{-} = g(\p_{+})$ with a gauge transformation of the form
$g = e^{2i\x}$ as discussed above, then we obtain the system
of differential equations
\begin{align}
    \label{eq:equivalent-connections}
	\del{\zbar} \log\phi_{-} & =  \del{\zbar} \log\phi_{+}\,-\,
	                               \del{\zbar}\,2i\x \\
	\label{eq:equivalent-fields}
	\del{z} \log\phi_{-}     & =  e^{2i\x}\,\del{z} \log\phi_{+}.
\end{align}
Note that equation (\ref{eq:equivalent-connections}) implies that
$\x$ is harmonic,
\begin{displaymath}
	\del{z}\del{\zbar}\,\x = 0,
\end{displaymath}
as it is the imaginary part of a holomorphic function.  Taking the 
conjugate of (\ref{eq:equivalent-fields}) and inserting it into 
(\ref{eq:equivalent-connections}), we can eliminate either
$\log\phi_{+}$ or $\log\phi_{-}$ from these equations to obtain
\begin{align*}
    \del{\zbar}\log\phi_{+} = 
        \frac{e^{2i\x}}{e^{2i\x}-1}\,\del{\zbar}\,2i\x
    &  \qquad
    \del{z}\log\phi_{+} =
       \frac{1}{e^{2i\x}-1}\,\del{z}\,2i\x \\
    \del{\zbar}\log\phi_{-} =
        \frac{1}{e^{2i\x}-1}\,\del{\zbar}\,2i\x
    &  \qquad
    \del{\zbar}\log\phi_{-} =
        \frac{e^{2i\x}}{e^{2i\x}-1}\,\del{z}\,2i\x.
\end{align*}   
These equations may also be written in either of the two simpler forms
\begin{align*}
    \del{\zbar}\log\phi_{\pm} = 
        \frac{e^{\pm i\x}}{\sin\x}\,\del{\zbar}\,\x
    &  \qquad
    \del{z}\log\phi_{\pm} =
       \frac{e^{\mp i\x}}{\sin\x}\,\del{z}\,\x
\end{align*}
or
\begin{align}
    \label{eq:dzbar}
    \del{\zbar}\log\phi_{\pm} & =
        \del{\zbar}\log\left( e^{\pm 2i\x}-1 \right) \qquad
    \del{z}\log\phi_{\pm}  = 
        \del{z}\log\left( e^{\mp 2i\x}-1 \right).
\end{align}
Substituting these formulae for the partial derivatives into
(\ref{eq:vortex-dbar}), we can
express the two gauge equivalent hyperbolic vortices completely
in terms of the gauge transformation without reference to their
super-potentials.
Furthermore, computing the Laplacian of the super-potentials
$\phi_{1},\phi_{2}$ in terms of the function $\x$, we have
\begin{equation*}
	\frac{1}{\phi_{\pm}}\,\del{\zbar}\del{z}\,\phi_{\pm}
	 =  
	    \del{\zbar}\log\phi_{\pm}\,\del{z}\log\phi_{\pm} \,+\,
	         \del{\zbar}\del{z}\log\phi_{\pm} 
	 =   
	    \frac{e^{\pm i\x}}{\sin\x}\,\del{z}\del{\zbar}\,\x.
\end{equation*}
Hence, the requirements that $\phi_{1},\phi_{2}$ be harmonic reduce
simply to the condition that $\x$ be harmonic, which we have already
established as a corollary to equation~(\ref{eq:equivalent-connections}).
We have thus proved

\begin{theorem} \label{vortex-construction}
    Let $\x$ be a real-valued harmonic function on the hyperbolic
    upper half-plane $\mathcal{H}^2$. The two pairs $(a_{+},\p_{+})$
    and $(a_{-},\p_{-})$ given by
    \begin{align*}
    	\dbar_{a_{\pm}} & =  \dbar\,+\,
    	                      \dbar\log\left( e^{\pm 2i\x}-1 \right)
    	                      \,+\, \frac{d\zbar}{z-\zbar} \\
    	\p_{\pm}        & =  i\,(z-\zbar)\,
    	                      \del{z}\log\left( e^{\mp 2i\x}-1 \right)
    \end{align*}
    then satisfy the hyperbolic vortex equations (\ref{eq:vortex-a})
    and (\ref{eq:vortex-b}) and are related by the gauge transformation
    $g = e^{2i\x}$. Conversely, any two gauge equivalent hyperbolic
    vortices constructed via the harmonic function ansatz of
    Theorem~\ref{theorem-vortex-ansatz} can be expressed in this form.
\end{theorem}


\subsection{The Unit Disc Model}

\label{unit-disc}

Until now, we have always used the upper half-plane model for hyperbolic 
space $\mathcal{H}^2$. In some circumstances, it will be more convenient to
use the unit disc model.  While the upper half-plane arises naturally by
the dimensional reduction technique discussed in the previous sections,
the calculations in the following section become much simpler and
exhibit significantly more symmetry if we can work on the unit disk.
Here we make the transition beween the two coordinate systems, showing 
how the formulae of the previous sections behave under the transformation.

Letting $z$ be the complex coordinate for the upper half-plane
and $w$ the coordinate for the unit disc, these two models are related
by the conformal transformation
\begin{equation} \label{eq:coordinate-transform}
	w = \frac{i-z}{i+z}, \qquad
	z = i\,\frac{1-w}{1+w}.
\end{equation}
This map takes the upper half-plane to the interior of the unit disc, 
mapping the real axis to the unit circle.  In particular, the point
$z = i$ maps to the origin, while the origin maps to $w = 1$ and the
point at infinity maps to $w = -1$.  The positive imaginary axis in 
$z$-coordinates maps to the interval $(-1,1)$ on the real axis in 
$w$-coordinates.

The hyperbolic metric on the unit disk is given by
\begin{displaymath}
	h = \frac{4}{\left(1-|w|^2\right)^2}\,dw\,d\wbar.
\end{displaymath}
The vortex equations remain fixed under this change of coordinates;
equation (\ref{eq:vortex-a}) is preserved because the transformation
is holomorphic, while equation (\ref{eq:vortex-b}) is a relationship
between coordinate-invariant scalar quantities. To compute the new
connection and Higgs fields generated by the harmonic function ansatz,
we first note that
\begin{displaymath}
	z - \zbar = 2i\,\frac{1-|w|^2}{|1+w|^2},
\end{displaymath}
and that the partial derivatives transform according to
\begin{displaymath}
	\del{w}     = -\frac{2i}{(1+w)^2}\,\del{z}, \qquad
	\del{\wbar} =  \frac{2i}{(1+\wbar)^2}\,\del{\zbar}.
\end{displaymath}
Converting the formula (\ref{eq:vortex-dbar})
for the hyperbolic vortex associated to a harmonic super-potential $\phi$
to $w$-coordinates on the unit disc, Theorem~\ref{theorem-vortex-ansatz}
then becomes

\begin{theorem} \label{w-theorem-vortex-ansatz}
    Given a positive real-valued super-potential $\phi$ over the hyperbolic
    disc $\mathcal{H}^{2}$, the $\U(1)$ connection and Higgs field pair 
    $(a, \p)$ defined by
    \begin{align} \label{eq:w-dbar}
	    \dbar_{a} = \dbar \,+\, \dbar\log\phi \,+\,
	                  \frac{d\wbar}{1-|w|^2}\,\frac{1+w}{1+\wbar} \qquad
	    \p        =-i\,(1-|w|^2)\,\frac{1+w}{1+\wbar}\,\del{w}\log\phi,
	\end{align}
	satisfies the vortex equations (\ref{eq:vortex-a}) and (\ref{eq:vortex-b})
	if and only if the super-potential $\phi$ is harmonic.
\end{theorem}

The Chern class $c_{1}(a)$ of this connection is given in these 
coordinates by
\begin{equation}\begin{split}
	c_{1}(a) =\frac{1}{2\pi}\int_{\mathcal{H}^{2}}
	               \ast_{h}\left(1-|\p|^{2}\right)  
	         =\frac{1}{2\pi}\int_{D^{2}}\left(
	               \frac{4}{\left(1-|w|^{2}\right)^{2}} -
	               4\,\left|\del{w}\log\phi\right|^{2}
	           \right) d\mu, \label{eq:c1-w}
\end{split}\end{equation}
where $d\mu$ is the volume element on the disc.


We now return to the simple hyperbolic vortex constructed
in \S\ref{vortex-tHooft}, which we obtained by a dimensional
reduction of the basic $c_{2}=1$ instanton. In $w$-coordinates on the
unit disc, the super-potential (\ref{eq:simple-potential}) becomes
\begin{equation}  \label{eq:w-simple-potential}
	\phi = 2\,\frac{1-|w|^{4}}{\left|1-w^{2}\right|^{2}}.
\end{equation}
Computing the partial derivatives of its logarithm, we obtain
\begin{equation*}
	\del{w}\log\phi  =  \frac{2w}{1-|w|^{4}}\,
	                      \frac{1-\wbar^{2}}{1-w^{2}} \qquad
	\del{\wbar}\log\phi  =  \frac{2\wbar}{1-|w|^{4}}\,
	                      \frac{1-w^{2}}{1-\wbar^{2}},
\end{equation*}
which when inserted into formula (\ref{eq:w-dbar}) above give the vortex
\begin{equation*}
    \dbar_{a}  =  \dbar + \frac{d\wbar}{1+|w|^{2}}\,
                            \frac{1+w}{1-\wbar}\, \qquad
	\p         =  -\frac{2iw}{1+|w|^{2}}\,\frac{1-\wbar}{1-w}.
\end{equation*}
Note that here the Higgs field $\p$ vanishes only at the origin,
where it has a simple zero, and so the Chern class of this vortex
should be $c_{1}(a) = 1$.

Using (\ref{eq:c1-w}) to explicitly calculate this Chern class,
we obtain the integral
\begin{displaymath}
	c_{1}(a) = \frac{1}{2\pi} \int_{D^{2}} \left(
	           \frac{4}{\left(1-|w|^{2}\right)^{2}} -
	           \frac{16\,|w|^{2}}{\left(1-|w|^{4}\right)^{2}}
                                           \right) d\mu.
\end{displaymath}
Taking polar coordinate $r,\theta$ on the unit disc, this
integral becomes
\begin{equation*}
	c_{1}(a)  =  \int_{0}^{1} \left(
	                   \frac{4\,r}{\left(1-r^2\right)^2} -
	                   \frac{16\,r^3}{\left(1-r^4\right)^2}
	                            \right) dr
	          =  \left[ \frac{2}{1-r^2} - \frac{4}{1-r^4}
	               \right]_{r=0}^{r=1}.
\end{equation*}
To evaluate this expression at $r=1$, we substitute $r=1-x$
and expand it about $x=0$, giving us
\begin{align*}
	\frac{2}{1-(1-x)^{2}} - \frac{4}{1-(1-x)^4}
	&= \frac{2}{2x-x^2+O(x^3)} - \frac{4}{4x-6x^2+O(x^3)} \\*
	&= \left(\frac{1}{x} + \frac{1}{2} + O(x)\right) -
	    \left(\frac{1}{x} + \frac{3}{2} + O(x)\right) \\*
	&= -1 + O(x) \rule{0in}{3.5ex}.
\end{align*}
We therefore see that the Chern class of this vortex is
indeed $c_{1}(a) = 1$ as we predicted by counting the zeros
of the Higgs field.

The formulae from \S\ref{gauge-transformation} giving the 
two super-potentials in terms of the gauge transformation remain 
unchanged, except for replacing all the $z$'s with $w$'s. In
particular, if the vortices determined by the super-potentials
$\phi_{+}$ and $\phi_{-}$ are gauge equivalent by a transformation 
of the form $g = e^{2i\x}$, then the differential equations
(\ref{eq:equivalent-connections}) and (\ref{eq:equivalent-fields}) become
\begin{align}
    \label{eq:w-equivalent-connections}
	\del{\wbar} \log\phi_{-} & =  \del{\wbar} \log\phi_{+}\,-\,
	                               \del{\wbar}\,2i\x \\
	\label{eq:w-equivalent-fields}
	\del{w} \log\phi_{-}     & =  e^{2i\x}\,\del{w} \log\phi_{+}.
\end{align}
The unit disc version of Theorem \ref{vortex-construction} is then
 
\begin{theorem} \label{w-vortex-construction}
    Let $\x$ be a real-valued harmonic function on the hyperbolic
    unit disc $\mathcal{H}^2$. The two pairs $(a_{+},\p_{+})$
    and $(a_{-},\p_{-})$ given by
    \begin{align*}
    	\dbar_{a_{\pm}} & =  \dbar\,+\,
    	                      \dbar\log\left( e^{\pm 2i\x}-1 \right)
    	                      \,+\,
    	                      \frac{d\wbar}{1-|w|^2}\,\frac{1+w}{1+\wbar} \\
    	\p_{\pm}        & =  -i\,(1-|w|^2)\,\frac{1+w}{1+\wbar}\,
    	                      \del{w}\log\left( e^{\mp 2i\x}-1 \right)
    \end{align*}
    then satisfy the hyperbolic vortex equations (\ref{eq:vortex-a})
    and (\ref{eq:vortex-b}) and are related by the gauge transformation
    $g = e^{2i\x}$. Conversely, any two gauge equivalent hyperbolic
    vortices constructed via the harmonic function ansatz of
    Theorem~\ref{w-theorem-vortex-ansatz} can be expressed in this form.
\end{theorem}

\section{Holonomy Singularity}

\subsection{The Forg\'{a}cs, Horv\'{a}th, Palla Instanton}

\label{fhp}

In this section, we construct the singular instanton described by 
P. Forg\'{a}cs, Z. Horv\'{a}th, and L. Palla in \cite{FHP1}.
In order to obtain a connection on $S^4\setminus S^2$ with a holonomy
singularity, Forg\'{a}cs {\em et al.\/} patch together two non-singular
connections on overlapping simply connected regions using a gauge
transformation. (This process is not unlike the clutching construction,
which creates ``twisted'' vector bundles given their local trivializations
and transition functions.) These two non-singular solutions are generated
by the harmonic function ansatz of Section~\ref{instanton-ansatz},
and since the super-potentials they use are $\SO(3)$-invariant,
the resulting connection can be analyzed in terms of the dimensional
reduction to hyperbolic space $\mathcal{H}^2$ discussed in Section 2.

Using the quaternionic notation $x = t + rQ$ with $t, r$ real, $r > 0$
and $Q$ pure imaginary satisfying $Q^{2} = -1$, we want to construct an
$\SO(3)$-invariant self-dual connection singular along the 2-sphere
$t=0, r=1$.  By dimensional reduction, this translates into a vortex
over the hyperbolic upper half-plane with non-trivial holonomy around
the point $z = i$.  If instead we work using the unit disc model of
hyperbolic space, our task takes the more symmetric form of finding a
hyperbolic vortex on the punctured disc with a holonomy singularity at
the origin.

We therefore set out to construct two gauge equivalent hyperbolic
vortices on the punctured disc, using the harmonic function ansatz
of \S\ref{vortex-ansatz}.  For the two simply connected
regions, we let $P_{1}$ be the disc with a cut along the positive
real axis, and let $P_{2}$ be the disc with a cut along the negative
real axis.  The areas of overlap are then the upper and lower half-discs,
excluding the real axis.  Let $(a_{1},\p_{1})$ and $(a_{2},\p_{2})$
be the hyperbolic vortices corresponding to the super-potentials
$\phi_{1}$ and $\phi_{2}$ on the regions $P_{1}$ and $P_{2}$ respectively.

In light of Theorem~\ref{vortex-construction}, we begin by examining
the gauge transformation between these two vortices, rather than focusing
on the super-potentials.  Here our gauge transformation $g$ is specified
by a real-valued harmonic function $\x$, and we take
\begin{displaymath}
	g = \left\{ \begin{array}{ll}
	                e^{+2i\x}  & \mbox{for $\Im w > 0$} \\
	                e^{-2i\x} & \mbox{for $\Im w < 0$}.
	            \end{array}
	    \right.
\end{displaymath}
In other words, on the lower half-disc we use the inverse of the gauge
transformation that we use on the upper half-disc.  In the notation
of \S\ref{gauge-transformation}, switching the sign of $2i\x$
simply interchanges the resulting gauge equivalent super-potentials
$\phi_{+}$ and $\phi_{-}$ determined by equation (\ref{eq:dzbar}).  Note that although the resulting $g$ is undefined
along the real axis, this does not pose a problem for our construction.
Indeed, we use $g$ directly only on regions excluding the real axis,
and we will find that the two hyperbolic vortices $(a_{1},\p_{1})$
and $(a_{2},\p_{2})$ constructed from $g$ are continuous across the
negative and positive axes respectively.

In \cite{FHP1}, Forg\'{a}cs {\em et al.\/} use the gauge transformation
specified by
\begin{equation} \label{eq:imlog}
    2\x = \left(\,
             \frac{\pi}{2}\,+\,2 \arctan \frac{T_{2}}{1-T_{1}}
                          \,+\,2 \arctan \frac{T_{1}}{1-T_{2}}
          \,\right).
\end{equation}
We will define the $T_{i}$ below, but before doing so we first
study the behavior of this gauge transformation for general values
of $T_{i}$.  In particular, we would like to coerce $T_{1}$ and
$T_{2}$ into being the real and imaginary parts of a holomorphic
(or anti-holomorphic) function $f(w)$.  Then where $|f(w)| = 1$,
the arguments of the arctans resemble the half-angle formula for
$\tan(\theta)$.  In such circumstances we have
\begin{equation*}
	\arctan \frac{T_{2}}{1-T_{1}} = \Im\log\,(1-\bar{f}) \qquad
	\arctan \frac{T_{1}}{1-T_{2}} = \Im\log\,(1+if),
\end{equation*}
and we can write $\x$ in terms of $f$ using
\begin{displaymath}
    2\x = \Im\log\left( i\,\frac{(1+if)^{2}}{(1-f)^{2}}
                 \right).
\end{displaymath}
We then see that $\x$ is indeed harmonic as it is the imaginary part
of a holomorphic (or anti-holomorphic) function.  Exponentiating, we
obtain
\begin{equation} \label{eq:exp-2ichi}
	e^{2i\x} = i\,\frac{(1-\bar{f})\,(1+if)}{(1-f)\,(1-i\bar{f})},
\end{equation}
and the expressions $e^{2i\x}-1$ and $e^{-2i\x}-1$ take the form
\begin{equation*}
	e^{+2i\x} - 1 = -\frac{(1-i)\,(1-f\bar{f})}{(1-f)\,(1-i\bar{f})}, \qquad
	e^{-2i\x} - 1 = -\frac{(1+i)\,(1-f\bar{f})}{(1-\bar{f})\,(1+if)},
\end{equation*}
which we shall use in equation (\ref{eq:dzbar}).
Note that if $|f| = 1$ then $\bar{f} = f^{-1}$, and we observe that
$e^{\pm 2i\x} = 1$.


In equation (\ref{eq:imlog}) above, the $T_{i}$ are defined by
\begin{align*}
	T_{1}&=\frac{1}{2\,S^{5}}\,\sqrt{(S+S_{-})^{2}-4}\,\left\{
	            \frac{1}{4} \left[4 - (S-S_{-})^2 \right]^{2}\,+\,
	            S^2 S_{-}^2\,-\,3\,(z+\zbar)^{2}
	                                   \right\} \\
	T_{2}&=\frac{1}{2\,S^{5}}\,\sqrt{4-(S-S_{-})^{2}}\,\left\{
	            \frac{1}{4} \left[4 - (S+S_{-})^2 \right]^{2}\,+\,
	            S^2 S_{-}^2\,-\,3\,(z+\zbar)^{2}
	                                   \right\}   
\end{align*}
with
\begin{displaymath}
	S = |z + i|, \qquad S_{-} = |z - i|,
\end{displaymath}
using the complex coordinate $z$ on the upper half-plane. With formulae
such as these, it is not surprising that the mathematical community was
incredulous.  Changing to the complex coordinate $w$ on the unit disc
via the conformal transformation (\ref{eq:coordinate-transform}), we have
\begin{displaymath}
	S     = \frac{2}{|w+1|}, \qquad
	S_{-} = 2\,\frac{|w|}{|w+1|} = |w| S.
\end{displaymath}
Calculating the various components of the $T_{i}$, we obtain
\begin{align*}
	\sqrt{(S+S_{-})^{2}-4} & = 
	    2\,\frac{|\sqrt{w}-\sqrt{\wbar}|}{|w+1|} =
	    \left|\sqrt{w}-\sqrt{\wbar}\right|S \\
	\sqrt{4-(S-S_{-})^{2}} & = 
	    2\,\frac{|\sqrt{w}+\sqrt{\wbar}|}{|w+1|} =
	    \left|\sqrt{w}+\sqrt{\wbar}\right|S \\
\end{align*}
and
\begin{displaymath}
	z + \zbar = -2i\,\frac{w-\wbar}{|w+1|^{2}}
	          = -\frac{i}{2}\,(w-\wbar)\,S^{2}.
\end{displaymath}
Putting these pieces together, the $T_{i}$ are given much more simply by
\begin{align*}
    T_{1} & =  \half \left|\sqrt{w}-\sqrt{\wbar}\right|
                \left[\,
                   \frac{1}{4}\left(\sqrt{w}+\sqrt{\wbar}\right)^{4} \,+\,
                   w\wbar \,+\, \frac{3}{4}\,(w-\wbar)^{2}
                \,\right] \\*
          & =  \half \left|w^{1/2}-\wbar^{1/2}\right|
                \left(
                    w^2 + w^{3/2}\wbar^{1/2} + w\wbar +
                          w^{1/2}\wbar^{3/2} + \wbar^{2}
                \right) \\*
          & =  \mp\frac{i}{2} \left(w^{5/2}-\wbar^{5/2}\right) 
          =  \left(\Im w^{5/2}\right)
                \left(\mbox{sign}\:\Im w^{1/2}\right),
\end{align*}
and
\begin{align*}
    T_{2} & =  \half \left|\sqrt{w}+\sqrt{\wbar}\right|
                \left[\,
                   \frac{1}{4}\left(\sqrt{w}-\sqrt{\wbar}\right)^{4} \,+\,
                   w\wbar \,+\, \frac{3}{4}\,(w-\wbar)^{2}
                \,\right] \\*
          & =  \half \left|w^{1/2}+\wbar^{1/2}\right|
                \left(
                    w^2 - w^{3/2}\wbar^{1/2} + w\wbar -
                          w^{1/2}\wbar^{3/2} + \wbar^{2}
                \right) \\*
          & =  \pm\half \left(w^{5/2}+\wbar^{5/2}\right) 
           =  \left(\Re w^{5/2}\right)
                 \left(\mbox{sign}\:\Re w^{1/2}\right).
\end{align*}
Hence, we see that $T_{1}$ and $T_{2}$ are indeed the imaginary and real
parts of the function $f(w)$ defined by
\begin{displaymath}
	f(w) = \left\{ \begin{array}{rcrcl}
	        w^{5/2}      && 0     & \le\;\;\arg w\;\;\le & \pi  \\
	        -\wbar^{5/2} && \pi   & \le\;\;\arg w\;\;\le & 2\pi
	               \end{array}
	       \right.
\end{displaymath}
on the region $P_{1}$ or equivalently
\begin{displaymath}
	f(w) = \left\{ \begin{array}{rcrcl}
	        \wbar^{5/2}  && -\pi  & \le\;\;\arg w\;\;\le & 0    \\
	        w^{5/2}      && 0     & \le\;\;\arg w\;\;\le & \pi  \\
	               \end{array}
	       \right.
\end{displaymath}
on the region $P_{2}$. This function $f(w)$ is holomorphic on the upper 
half-disc and anti-holomorphic on the lower half-disc, and we note that
$f(w)$ is well defined and continuous over the whole unit disc.

We now have all that we need to calculate the gauge equivalent
connections and Higgs fields $a_{i},\p_{i}$ satisfying
$a_{2} = g(a_{1})$ and $\p_{2} = g(\p_{1})$.
We first consider
the super-potential $\phi_{1}$ defined on the region
$0 < \arg w < 2\pi$.  Using the notation of
\S\ref{gauge-transformation}, on the upper half-disc
we have $\phi_{1} = \phi_{+}$ and $f = w^{5/2}$, giving us
\begin{displaymath}
	\del{\wbar}\log\phi_{1}
	= \del{\wbar}\log\frac{1-|w|^{5}}
	                      {(1-w^{5/2})\,(1-i\,\wbar^{5/2})}
	=\frac{5}{2}\,\frac{i\,\wbar^{3/2}}{1-|w|^5}\,\frac{1+i\,w^{5/2}}{1-i\,\wbar^{5/2}}.
\end{displaymath}
On the lower half-disc, we obtain the same expression
\begin{displaymath}
	\del{\wbar}\log\phi_{1}
	= \del{\wbar}\log\frac{1-|w|^{5}}
	                      {(1+w^{5/2})\,(1-i\,\wbar^{5/2})}
	= \frac{5}{2}\,\frac{i\,\wbar^{3/2}}{1-|w|^5}\,\frac{1+i\,w^{5/2}}{1-i\,\wbar^{5/2}},
\end{displaymath}
although this time we take $\phi_{1}=\phi_{-}$ and $f=-\wbar^{5/2}$.
Similarly, considering the super-potential $\phi_{2}$ defined
for $-\pi < \arg w < \pi$, on the upper half-disc with
$\phi_{2} = \phi_{-}$ and $f = w^{5/2}$ we have
\begin{displaymath}
	\del{\wbar}\log\phi_{2}
	= \del{\wbar}\log\frac{1-|w|^{5/2}}
	                      {(1-\wbar^{5/2})\,(1+i\,w^{5/2})}
	= \frac{5}{2}\,\frac{\wbar^{3/2}}{1-|w|^5}\,\frac{1-w^{5/2}}{1-\wbar^{5/2}},
\end{displaymath}
while on the lower half-disc, we again get the same expression
\begin{displaymath}
	\del{\wbar}\log\phi_{2}
	= \del{\wbar}\log\frac{1-|w|^{5/2}}
	                      {(1-\wbar^{5/2})\,(1-i\,w^{5/2})}
	= \frac{5}{2}\,\frac{\wbar^{3/2}}{1-|w|^5}\,\frac{1-w^{5/2}}{1-\wbar^{5/2}},
\end{displaymath}
taking $\phi_{2} = \phi_{+}$ and $f = \wbar^{5/2}$. Hence, as we
predicted, the vortices determined by $\phi_{1}$ and $\phi_{2}$ extend
continuously across the negative and positive real axes respectively,
even though the gauge transformation $g$ between them does not.

For the purposes of our vortex construction, we need not know the 
super-potentials explicitly; rather, all we require are the complex partial
derivatives of their logarithms which we computed above. Nevertheless,
here we present the super-potentials as given in \cite{FHP1}. There the two
super-potentials take center stage, defined on their own instead of
being constructed from the gauge transformation as we have done.
Forg\'{a}cs {\em et al.\/} define $\phi_{1}$, $\phi_{2}$ by%
\SSfootnote{Actually, \cite{FHP1} uses $+2\,S^{5}\,T_{1}$ and $+2\,S^{5}\,T_{2}$
          in the denominators of $\phi_{1}$ and $\phi_{2}$ respectively.
          The resulting super-potentials are still harmonic, and they do
          indeed generate gauge equivalent hyperbolic vortices. However,
          this choice of sign is not consistent with the gauge
          transformation they use, given here by (\ref{eq:imlog}).}
\begin{displaymath}
	\phi_{1}=\frac{S^{5}-S_{-}^{5}}{S^{5}+S_{-}^{5}-2\,S^{5}\,T_{1}},
	    \qquad
	\phi_{2}=\frac{S^{5}-S_{-}^{5}}{S^{5}+S_{-}^{5}-2\,S^{5}\,T_{2}}.
\end{displaymath}
Simplifying these expressions and writing them using the unit disc model 
of hyperbolic space, we obtain
\begin{align} \label{eq:w-phi1}
	\phi_{1}&=\frac{1-|w|^{5}}{1\,-\,2\,\Im w^{5/2}\,+\,|w|^5}
	  =\frac{1-|w|^{5}}{\left|1+i\,w^{5/2}\right|^{2}} \\
	             \label{eq:w-phi2}
	\phi_{2}&=\frac{1-|w|^{5}}{1\,-\,2\,\Re w^{5/2}\,+\,|w|^5}
	  =\frac{1-|w|^{5}}{\left|1-w^{5/2}\right|^{2}}.
\end{align}
The reader may want to verify that the partial derivatives
$\del{\wbar}\log\phi_{1}$ and $\del{\wbar}\log\phi_{2}$
agree with those calculated on the previous page.

We now compute the Chern class $c_{1}$ of the vortex patched together from
the two super-potentials $\phi_{1}$ and $\phi_{2}$. Using equation
(\ref{eq:c1-w}), we have
\begin{displaymath}
	c_{1}(a) = \frac{1}{2\pi}\int_{D^{2}}\left(
	               \frac{4}{(1-|w|^{2})^{2}} -
	               \frac{25\,|w|^{3}}{(1-|w|^{5})^{2}}
	                                     \right) d\mu.
\end{displaymath} 
Taking polar coordinates $r,\theta$ on the disc $D^{2}$, this integral 
becomes
\begin{equation*}
	c_{1}(a)=\int_{0}^{1}\left(
	               \frac{4\,r}{(1-r^{2})^{2}} -
	               \frac{25\,r^{4}}{(1-r^{5})^{2}}
	                         \right) dr 
	        =\left[\frac{2}{1-r^{2}} - \frac{5}{1-r^{5}}
	           \right]_{r=0}^{r=1},
\end{equation*}
and evaluating this expression by the method used at the end of
\S\ref{unit-disc} shows that $c_{1}(a) = 3/2$. We can also arrive
at this same result by naively counting the zeros (with multiplicity)
of the Higgs field, as both $\del{w}\log\phi_{1}$ and
$\del{w}\log\phi_{2}$ vanish to order $3/2$ at the origin
but are otherwise nonzero.


\subsection{A Family of Singular Vortices}
In this section we will generalize the construction of \cite{FHP1} to
produce a family of hyperbolic vortices with varying Chern class
$c_{1}$. This family includes both the standard $c_{1} = 1$
vortex of \S\ref{unit-disc} and the fractionally charged
vortex of the previous section, as well as vortices with arbitrary
real $c_{1}$. We will continue to work using the unit disc model of
hyperbolic space.

\newcommand{\e}{\epsilon}
\newcommand{\ebar}{\bar{\epsilon}}


Observing the resemblance between the two super-potentials
(\ref{eq:w-simple-potential}) and (\ref{eq:w-phi2}), we consider
a more general super-potential of the form
\begin{equation} \label{eq:phi-c}
	\phi = \frac{1-|w|^{2c}}{\left|1-w^{c}\right|^{2}},
\end{equation}
where $c$ is a nonzero real constant. For non-integral $c$, this
$\phi$ is not well defined over the whole unit disc; rather, we
must restrict it to a simply-connected cut disc. If we loop around
the origin once in the positive direction, crossing our cut in
the unit disc, then this super-potential becomes
\begin{equation} \label{eq:phi-c-epsilon}
	\phi' = \frac{1-|w|^{2c}}{\left|1-\e w^{c}\right|^{2}},
\end{equation}
introducing a factor of $\e = e^{2\pi ic}$ in the denominator.
In both of these cases, we note that $\phi$ and $\phi'$ vanish on
on the unit circle, except at the roots of $w^{a} = 1$
(or $w^{a} = \ebar$) where they have simple poles.


Taking the logarithmic derivatives of the super-potential $\phi$,
we obtain
\begin{equation*}
	\del{w}\log\phi      =  \frac{c\,w^{c-1}}{1-|w|^{2c}}\,
	                          \frac{1-\wbar^{c}}{1-w^{c}} \qquad
	\del{\wbar}\log\phi  =  \frac{c\,\wbar^{c-1}}{1-|w|^{2c}}\,
	                          \frac{1-w^{c}}{1-\wbar^{c}},
\end{equation*}
which we can use in equation (\ref{eq:w-dbar}) 
to construct a hyperbolic vortex $(a,\p)$. After looping around the
origin, the logarithmic derivatives become
\begin{equation*}
	\del{w}\log\phi'      =  \frac{c\,\e\,w^{c-1}}{1-|w|^{2c}}\,
	                           \frac{1-\ebar\,\wbar^{c}}{1-\e\,w^{c}} \qquad
	\del{\wbar}\log\phi'  =  \frac{c\,\ebar\,\wbar^{c-1}}{1-|w|^{2c}}\,
	                           \frac{1-\e\,w^{c}}{1-\ebar\,\wbar^{c}},
\end{equation*}
and we let $(a',\p')$ be the corresponding hyperbolic vortex.


In order to construct a single hyperbolic vortex over the whole of
the unit disc, we would like to find a gauge transformation $g$
taking the vortex $(a,\p)$ to the vortex $(a',\p')$. Using equation 
(\ref{eq:w-equivalent-fields}) for gauge equivalent Higgs fields,
we have
\begin{displaymath}
	\del{w}\log\phi' = g\,\del{w}\log\phi,
\end{displaymath}
and so our gauge transformation must be
\begin{equation} \label{eq:g}
	g = \e\,\frac{1-\ebar\,\wbar^{c}}{1-\e\,w^{c}}\,
	        \frac{1-w^{c}}{1-\wbar^{c}}.
\end{equation}
We must also verify that this gauge transformation $g$ satisfies equation
(\ref{eq:w-equivalent-connections}) for gauge equivalent connections
\begin{displaymath}
	\del{\wbar}\log\phi' = \del{\wbar}\log\phi - \del{\wbar}\log g,
\end{displaymath}
which we leave to the reader. Hence when we loop around the origin, we 
obtain a vortex that is gauge equivalent to our original one, and so
this vortex is well defined over the punctured disc.

Suppose that instead of defining $\epsilon = e^{2\pi ic}$, we use
an arbitrary constant $\epsilon$ with $|\epsilon| = 1$ in equations
(\ref{eq:phi-c-epsilon}) and (\ref{eq:g}). In this case, the gauge
transformation $g$ still maps between the two vortices corresponding
to the super-potentials $\phi$ and $\phi'$. Indeed, if we take
$c = 5/2$ and $\epsilon = -i$, then the resulting $\phi$, $\phi'$,
and $g$ correspond to the super-potentials (\ref{eq:w-phi2}) and
(\ref{eq:w-phi1}) and the gauge transformation (\ref{eq:exp-2ichi})
we used in \S\ref{fhp}.


Writing out the connection $a$ and Higgs field $\p$ of this vortex
explicitly using formula (\ref{eq:w-dbar}),
we obtain
\begin{align*}
	\dbar_{a}&= \dbar + \left(
	             \frac{c\,\wbar^{c-1}}{1-|w|^{2c}}\,
	             \frac{1-w^{c}}{1-\wbar^{c}} +
	             \frac{1}{1-|w|^{2}}\,\frac{1+w}{1+\wbar}
	                    \right) d\wbar \\
	\p       &= -ic\,w^{c-1}\,\frac{1-|w|^{2}}{1-|w|^{2c}}\,
	             \frac{1+w}{1+\wbar}\,
	             \frac{1-\wbar^{c}}{1-w^{c}}.
\end{align*}
The Chern class $c_{1}(a)$ of this vortex is given by equation
(\ref{eq:c1-w}), yielding the integral
\begin{align*}
	c_{1}(a) & =  \frac{1}{2\pi}\int_{D^{2}}\left(
	               \frac{4}{\left(1-|w|^2\right)^{2}} -
	               \frac{4c^{2}|w|^{2c-2}}{\left(1-|w|^{2c}\right)^{2}}
	                                         \right) d\mu \\
	         & =  \int_{0}^{1}\left(
	               \frac{4r}{(1-r^{2})^{2}} -
	               \frac{4c^{2}r^{2c-1}}{(1-r^{2c})^{2}}
	                           \right) dr
	          =  \left[
	               \frac{2}{1-r^{2}} - \frac{2c}{1-r^{2c}}
	               \right]_{r=0}^{r=1},
\end{align*}
where we have used polar coordinates $r,\theta$ on the unit disc. 
Evaluating this final expression by the method used at the end of
\S\ref{unit-disc}, we have $c_{1}(a) = c - 1$.


For $c = 2$, this construction yields the standard $c_{1} = 1$
hyperbolic vortex which we discussed in \S\ref{unit-disc}.
If we take $c = 5/2$ then we obtain the $c_{1} = 3/2$ vortex
given by the super-potential (\ref{eq:w-phi2}) on the cut
disc $P_{2}$. Note that with our construction, it is no longer
necessary to cover the disc with two overlapping regions as
we did in \S\ref{fhp}. Rather, it is sufficient to take
a single vortex on the cut disc and then study how it behaves
across that cut. For integer values of $c$, the vortex is
continuous across the cut, while for other values the vortex
changes by a gauge transformation. The flat $c=1$ vortex in
this family is
\begin{equation*}
	\dbar_{a} = \dbar + \frac{2\,d\wbar}{1-\wbar^{2}} \qquad
	\p        = -i\,\frac{1+w}{1+\wbar}\,\frac{1-\wbar}{1-w},
\end{equation*}
which we readily see has $|\p| = 1$ and $F_{a} = 0$, and we
therefore have $c_{1} = 0$ as expected.  We observe that this
vortex is equivalent to the standard flat vortex $(a=0, \p=1)$
using the gauge transformation (\ref{eq:g}) with $\epsilon = -1$.

To compute the holonomy around loops circling the origin, we
introduce polar coordinates $r, \theta$ on the unit disc. In
these coordinates, the complex differentials $dw$ and $d\wbar$
are
\begin{displaymath}
	dw     = \frac{w}{|w|}\,dr + iw\,d\theta, \qquad
	d\wbar = \frac{\wbar}{|\wbar|} \, dr - i\wbar\,d\theta.
\end{displaymath}
In a small neighborhood of the origin, our connection is approximated by
\begin{equation*}
	a \approx c \left( \wbar^{c-1}\,d\wbar - w^{c-1}\,dw \right) 
	   =      -ic \left( w^{c} + \wbar^{c} \right) d\theta + \cdots
\end{equation*}
(only the $d\theta$ term is needed for calculating the holonomy).
The contribution to the holonomy around the circle $|w| = r$ due to
the connection $a$ is then given by the loop integral
\begin{align*}
    \oint_{|w|=r}a &= \int_{0}^{2\pi}\!-ic r^{c}
	                   \left( e^{ic\theta} + e^{-ic\theta}
	                   \right) d\theta \\
	   &=r^{c} \left( -e^{2\pi ic} + e^{-2\pi ic} \right) 
	   =-2 r^{c} \sin 2\pi ic.
\end{align*}
For $c > 0$, this expression vanishes as $r \rightarrow 0$, and thus
the limit holonomy around the origin comes entirely from the gauge
transformation $g$ given by (\ref{eq:g}). Near the origin, we have
$g\approx\epsilon=e^{2\pi ic}$. Hence, the limit of the holonomy around
small loops centered at the origin is $e^{2\pi ic}$.

By the dimensional reduction technique of Section 2, our family of
singular hyperbolic vortices corresponds to a similar family of
$\SO(3)$-invariant self-dual connections over $S^{4}\setminus S^{2}$.
Using the results and quaternionic notation of \S\ref{gauge-transformation},
we see that the holonomy around small loops linking the singular surface
$S^{2}$ is $e^{\pi Qc}$. Viewing these connections as $\SU(2)$ connections
on a bundle $E$, we then obtain a splitting $E=L\oplus L^{\ast}$ on a 
neighborhood of $S^{2}$ with respect to which the holonomy takes the
standard form
\begin{displaymath}
	\left( \begin{array}{cc}
	    e^{2\pi i\alpha} & 0 \\
	    0                & e^{-2\pi i\alpha}
	\end{array} \right)
\end{displaymath}
for a constant $\alpha$ in the range $[0,1/2)$.  In the $\SO(3)$-symmetric
case, the complex line bundle $L$ on $S^{2}$ has Chern class $c_{1}(L)=-1$.
For our family of singular solutions, the holonomy parameter is
$\alpha=(c-\lfloor c\rfloor)/2$, where $\lfloor c\rfloor$ is the greatest
integer less than or equal to $c$.

\end{document}